\documentclass[11pt,letterpaper]{article}

\usepackage[utf8]{inputenc}
\usepackage[T1]{fontenc}
\usepackage{mathtools}
\usepackage{amsmath,amssymb, mathrsfs,amsthm}
\usepackage{hyperref}
\usepackage[a4paper, total={15cm, 22cm}]{geometry}
\usepackage{tikz-cd}
\usepackage{xcolor}
\usepackage{todonotes}

\DeclarePairedDelimiter{\norm}{\lVert}{\rVert}

\newtheorem{theorem}{Theorem}[section]
\newtheorem{lemma}[theorem]{Lemma}
\newtheorem{conjecture}[theorem]{Conjecture}
\newtheorem{proposition}[theorem]{Proposition}
\newtheorem{corollary}[theorem]{Corollary}
\newtheorem{definition}[theorem]{Definition}
\newtheorem{notation}[theorem]{Notation}
\newtheorem{remark}[theorem]{Remark}
\newtheorem{example}[theorem]{Example}

\newcommand{\bbP}{\mathbb{P}}
\newcommand{\bbR}{\mathbb{R}}
\newcommand{\bbN}{\mathbb{N}}
\newcommand{\bbQ}{\mathbb{Q}}
\newcommand{\bbC}{\mathbb{C}}
\newcommand{\bbZ}{\mathbb{Z}}

\newcommand{\cX}{\mathcal{X}}
\newcommand{\cY}{\mathcal{Y}}
\newcommand{\cO}{\mathcal{O}}
\newcommand{\cJ}{\mathcal{J}}

\newcommand{\frX}{\mathfrak{X}}
\newcommand{\frY}{\mathfrak{Y}}
\newcommand{\fra}{\mathfrak{a}}
\newcommand{\frm}{\mathfrak{m}}

\newcommand{\frj}{\mathfrak{j}}
\newcommand{\frU}{\mathfrak{U}}

\renewcommand{\div}{\ensuremath{\operatorname{div}}}

\newcommand{\Div}{\ensuremath{\operatorname{Div}}}
\newcommand{\CDiv}{\ensuremath{\operatorname{CDiv}}}

\newcommand{\Spec}{\ensuremath{\operatorname{Spec}}}
\newcommand{\ord}{\ensuremath{\operatorname{ord}}}
\newcommand{\mult}{\ensuremath{\operatorname{mult}}}
\newcommand{\vol}{\ensuremath{\operatorname{vol}}}

\newcommand{\red}{\ensuremath{\operatorname{red}}}
\newcommand{\cent}{\ensuremath{\operatorname{center}}}
\newcommand{\Env}{\ensuremath{\operatorname{Env}}}
\newcommand{\Mov}{\operatorname{Mov}}
\newcommand{\Nef}{\operatorname{Nef}}
\newcommand{\trdeg}{\operatorname{tr.deg.}}
\newcommand{\ratrk}{\operatorname{rat.rk.}}

\newcommand{\length}{\operatorname{length}}
\newcommand{\vsq}{\mathbf{v}^2}

\newcommand{\ordnorm}{\ensuremath{\theta}}

\title{\textbf{b}-divisorial valuations and Berkovich positivity functions}
\author{Joaquim Ro\'e\footnote{Partially supported by Grant PID2020-116542GB-I00 funded by MCIN/AEI/10.13039/501100011033.} {} and Stefano Urbinati\footnote{Partially supported by the project DM737 “Birational geometry, Fano manifolds and torus actions” CUP G25F21003390007, European Union funds, NextGenerationEU.\newline\phantom{xxx}2020 MSC Primary 14E99; Secondary 14C20, 14F18, 14B10, 12J20, 13A18.
} } 

\date{}

\begin{document}
\maketitle

\begin{abstract}
    We prove semicontinuity properties for local positivity invariants of big and nef divisors.
    The usual definition of Seshadri constant and asymptotic order of vanishing along a subvariety is extended to include all seminorms in the Berkovich space, and
    we obtain semicontinuity of such constants as a function of the center seminorm.
    
    We use Shokurov's language of \textbf{b}-divisors; to each seminorm there is an associated \textbf{b}-divisor which can be used to translate questions about positivity into questions about the shape of certain cones of \textbf{b}-divisors. 
    The theory works especially well for what we call \textbf{b}-divisorial valuations, a natural extension of the notion of divisorial valuations which encompasses e.g., all Abhyankar valuations.

\end{abstract}

\section{Introduction}

A profound, yet often elusive, connection links the Berkovich analytification, the Riemann-Zariski space, and the theory of \textbf{b}-divisors in birational geometry. 
This work, motivated by a recurring analytic pathology encountered in problems connected to Nagata's celebrated conjecture on singular plane curves, aims to render this connection explicit from a geometric and birational viewpoint, and to show its signification for the study of positivity invariants such as Seshadri constants. 
We explore how these global geometric invariants behave under variation in high-dimensional spaces of valuations, and how their discontinuities reflect the underlying structure of these spaces.

In the work of \cite{DHKRS17} a key function was constructed on a branch of the Berkovich space of the projective plane, attached to quasi-monomial valuations. 
The authors demonstrated how continuity of this function could be used to tackle Nagata's conjecture; on the other hand the function extends to the whole Berkovich space and behaves semicontinuously there. 
Surprisingly, an identical function emerged in the seemingly unrelated domain of symplectic geometry, specifically in the study of the \emph{ellipsoid packing} problem \cite{Dusa}, suggesting a subtle correspondence between algebraic and symplectic invariants.
A third intimately connected instance of an important geometric invariant naturally parametrized by such a valuation space occurs in the theory of Newton-Okounkov bodies. These are convex bodies associated to a line bundle on a projective variety endowed with a valuation of full rank, designed to encode asymptotic section data. In \cite{CFKLRS17}, a fundamental pathology was observed: while the volume of a Newton-Okounkov body remains invariant under a change of valuation, its \emph{shape} can change \emph{discontinuously} even as the valuation moves along a branch of the Berkovich space. This \emph{mutation} phenomenon is also observed in the symplectic setting \cite{MS25} and might suggest that the standard topology on the space of valuations is not fine enough to control the geometric data it is meant to parameterize.

These examples, while mysterious and currently largely confined to a 2-dimensional setting, point towards a fundamental connection between the topology of valuation spaces and positivity measures in birational geometry. This work aims to provide a unifying explanation for these pathological phenomena through the lens of classical valuation theory and birational geometry. We lay a geometric groundwork that clarifies the nature of these discontinuities, paving the way for future explorations, particularly into the connections with symplectic geometry.

Let \( X \) be a smooth projective variety over a closed field \( k \). We explore the relationship between two fundamental, non-Archimedean objects:
\begin{itemize}
    \item The \emph{Berkovich analytification} \( \frX^{\mathrm{an}} \), whose points correspond to multiplicative seminorms on the function field \( K(X) \) extending the base field norm, which we take to be the trivial norm. This space is endowed with a natural and well-behaved (e.g., locally path-connected) topology, making it a central object in non-Archimedean geometry.
    \item The \emph{Zariski-Riemann space} \( \cX \), defined as the space of all valuation rings on \( K(X) \) that are centered on \( X \), or equivalently, as the space of all equivalence classes of valuations \( \nu: K(X)^{\times} \to \Gamma \) on \( K(X) \) centered on \( X \).
\end{itemize}

We approach the space \( \cX \) from two complementary perspectives. The first is the classical, algebraic definition, which offers conceptual clarity but often obscures the underlying birational geometry. While this viewpoint yields important general results, such as the density of divisorial valuations, the non-noetherian geometric structure of \( \cX \) remains difficult to grasp. 

The second perspective, championed by Shokurov (cf. \cite{BdFF12}), provides a more geometric and birational characterization. Here, \( \cX \) is realized as the projective limit of all proper birational models of \( X \) over its base:
\[
\cX = \varprojlim_{X' \to X} X',
\]
where the limit is taken over all proper birational morphisms \( X' \to X \) with \( X' \) normal. This definition highlights \( \cX \) as a universal object in birational geometry. The bridge between these two definitions is furnished by the fact that the towers of birational maps are intrinsically governed by the inclusion relations of local rings and valuation rings.

The novel contribution of this work is to build a bridge between these points of view by systematically analyzing the \emph{trace} (or \emph{support}) of valuation ideals (relative to a base model $X$) on the projective system of models defining the Zariski-Riemann space \( \cX \). This leads naturally to the construction of a \textbf{b}-divisor $D_\xi$ ---a concept central to modern birational geometry, representing a collection of divisors on all birational models that are compatible under pushforward--- naturally associated to any seminorm $\xi$. 
Our work relies heavily on a structure introduced by Boucksom--De Fernex--Favre \cite{BdFF12} and Jonsson--Musta\c{t}\u{a} \cite{JM12}; indeed, $D_\xi$ is the \textbf{anti-effective} \textbf{b}-divisor associated, in the sense of \cite{BdFF12}, to the graded sequence of ideals of the valuation in the spirit of \cite{JM12} (see Section \ref{sec:d-xi} for details).
By elucidating to what extent the data of a valuation is encoded in this limiting divisor on all possible birational models, we provide a clear geometric framework for understanding the discontinuities observed in the motivating examples. This approach reveals the intricate interplay between analysis, algebra, and geometry in the study of valuation spaces, positioning the space of \textbf{b}-divisors as a well-suited geometric setting for analyzing the continuity of positivity invariants.

In our approach, the class of valuations for which the associated \textbf{b}-divisor allows one to recover the valuation is of great importance. 
We propose to call these valuations \emph{\textbf{b}-divisorial}, and they strictly include quasimonomial, or Abhyankar, valuations, and hence all divisorial valuations. 
Arguably, the excellent properties of Abhyankar valuations in birational geometry (see \cite{KK05}, \cite{ELS03}, \cite{JM12} among many) derive from their being \textbf{b}-divisorial.
We conjecture in \ref{con:characterize-b-divisorial} that \textbf{b}-divisorial valuations are characterized by the inequality $D_\xi\ne 0$. 
In the direction of characterizing \textbf{b}-divisorial valuations, we prove the following.

\begin{theorem}
 [\ref{thm:b-divisorial surfaces},\ref{cor:characterize-b-divisorial-surfaces}] \label{introthm:b-divisorial}
    Let $v$ be a valuation of $K(X)$ and $\xi$ the associated norm.
    If $\dim X=2$, then $v$ is \textbf{b}-divisorial if and only if $D_\xi\ne 0.$

    In general, if $v$ has sublinear log-discrepancy then $v$ is \textbf{b}-divisorial.
\end{theorem}

Valuations with sublinear log-discrepancy include all valuations with finite log-dis\-crep\-an\-cy \cite{JM12}, see Definition \ref{def:sublinear-discrepancy}.

Not long after this work was made available on arXiv, Lu Qi announced in \cite{qi26} a proof of our conjecture, i.e., that $v$ is \textbf{b}-divisorial if and only if $D_\xi\ne 0$ in arbitrary dimension, leveraging the \emph{saturation} property of the associated filtrations, see \cite{BLQ24}.

\medskip 

For our purposes it is crucial to understand how the \textbf{b}-divisor $D_\xi$ varies as a function of $\xi$. In this regard, we show the following.
\begin{theorem}  [\ref{thm:lowcont}, \ref{cor:delta-notation}] \label{introthm:D-xi-lowcont}
	The map $D:\frX \rightarrow \Div_\bbR(\cX)$ given by $\xi \mapsto D_\xi$ is lower semicontinuous.

    Moreover, for every model $X_\pi$ above $X$, every point $p\in X_\pi$, and every local system of parameters $z$, the map $D$ restricted to the set $\Delta$ of quasimonomial valuations with respect to $z$ is continuous in the interior of $\Delta$. 
\end{theorem}

Then, with these results at hand, in sections \ref{sec:positivity} and \ref{sec:surfaces} we address our original problem, which is to study \emph{local positivity invariants}, such as Seshadri constants and asymptotic orders of vanishing, as functions of the associated valuation. First we extend the classical definitions of these positivity invariants to general seminorms; if $D$ is an ample divisor on $X$, its Seshadri constant at a \textbf{b}-divisorial valuation $v_\xi$ is 
   $$\varepsilon(D, v_\xi):=\sup \left\{t \, |\, (D +t D_{\xi}) \text{ is nef }\right\},$$
and its asymptotic order of vanishing is
   $$\omega(D, v_\xi):=\sup \left\{t \, |\, (D +t D_{\xi}) \text{ is pseudoeffective }\right\}.$$
Finally, our main result, which sheds some light on the continuity and discontinuity phenomena that motivated us, is as follows.

\begin{theorem} [\ref{thm:waldshmidt-continuity}, \ref{thm:toric-seshadri-continuity}, \ref{thm:continuity-seshadri-surfaces}]
    \label{introthm:positivity-semicontinuity}
    Fix $D$ an ample divisor, and consider the functions $\frX\rightarrow \bbR\cup\{\infty\}$ 
    given by $\xi\mapsto\varepsilon(D,\xi)$ and $\xi\mapsto\omega(D,\xi)$.
    \begin{itemize}
        \item $\omega$ is lower semicontinuous, and it is continuous in the interior of every set $\Delta$ of quasimonomial valuations as above.
        \item $\varepsilon$ is continuous in the interior of every set $\Delta$ of quasimonomial valuations as above.
        \item If $X$ is a surface, and $p\in X$ a smooth point, then $\varepsilon$ is lower semicontinuous over the subset $\frX_p$ of valuations centered at $p$.
    \end{itemize}
\end{theorem}

Note that $\varepsilon$ is not globally lower (nor upper) semicontinuous, see Example \ref{exa:non-semicontinuous-seshadri}. 
Since both asymptotic order of vanishing and Seshadri constant are reflected in the shape of Newton--Okounkov bodies, our results are also significant for the understanding of the discontinuities observed in \cite{CFKLRS17}.
We expect to go back to the issue of semicontinuity for Newton--okounkov bodies in later work. 
The use of \textbf{b}-divisorial valuations should also help in this context by unifying the analysis of bodies determined by valuations of maximal rank which are not flag valuations, especially on surfaces.

Section \ref{sec:preliminaries} is devoted to fixing the language and notation for both types of valuation spaces we use, including some results on the theory of \textbf{b}-divisors from \cite{BdFF12} and on multiplicative ideals from \cite{JM12} which will be used extensively.
In section \ref{sec:d-xi} we introduce the divisor $D_\xi$ associated to a seminorm and we initiate its study, including the proof of Theorem \ref{introthm:b-divisorial}.
Section \ref{sec:continuity_Dxi} is devoted to the proof of Theorem \ref{introthm:D-xi-lowcont}.
In Section \ref{sec:positivity} we define and study the positivity invariants associated to seminorms, and prove most of Theorem \ref{introthm:positivity-semicontinuity}. 
Finally, in Section \ref{sec:surfaces} we deal with the case of surfaces, in which it is possible to give an explicit determination of $D_\xi$ in terms of the sequence of centers of a valuation, and we prove the part of Theorem \ref{introthm:positivity-semicontinuity} dealing with surfaces.

\section{Preliminaries on valuation spaces}
\label{sec:preliminaries}

Fix an algebraically closed field $K$ of characteristic zero, 
endowed with the trivial valuation.
Fix $X$ a smooth irreducible projective variety over $K$. 
We will be working on two spaces whose points may be identified with valuations of $K(X)$, namely the Berkovich analytification and the Riemann--Zariski space of $X$.  
We begin by recalling the basics and setting up notation for these spaces.

We follow the usual conventions of writing valuations additively and seminorms multiplicatively.

\subsection{Berkovich analytification}
\label{sec:intro-Berkovich}

Let $\frX$ stand for the Berkovich analytification of $X$. 
It is a Berkovich analytic space (see \cite{Ber90}, \cite{Ber93}, \cite{BCDKT08}, \cite{Tem15}, \cite{Thu07} for general definitions and results) 
whose points are identified with real multiplicative seminorms on $\cO_X$.
More precisely, a point $\xi\in\frX$ is given by an affine open set $U_\xi\subset X$ and a multiplicative seminorm $\norm{\cdot}_\xi:\cO_X(U_\xi)\rightarrow \bbR_{>0}$ whose restriction to $K$ is trivial, with the natural identifications on overlaps.

\textbf{The kernel map and the reduction map}. The analytification $\frX$ comes with a map $\ker:\frX\rightarrow X$ which sends each $\xi$ to the (schematic) point associated to its kernel, the prime ideal in $\Spec \cO_X(U_\xi)\subset X$ formed by those $f\in \cO_X(U_\xi)$ with $\norm{f}_\xi=0$.
The fiber of $\ker$ over a point $p\in X$ can be identified with the set of norms on $\kappa(p)$, the field of rational functions of the closed subvariety $\overline{p}$, which have trivial restriction to $K$.
Giving a norm $\norm{\cdot}$ on $\kappa(p)$ is equivalent to giving a real valuation $v$ on $\kappa(p)^\times$, related by $v(f)=-\log\norm{f}$.
Thus a point $\xi$ in $\frX$ can be completely identified with a pair $(V,v)=(V_\xi,v_\xi)$ where $V\subseteq X$ is a closed subvariety and $v:K(V)^\times \rightarrow \mathbb{R}$ is a valuation. Then $\ker(V_\xi,v_\xi)$ is simply the generic point of $V_\xi$, and
$$\norm{f}_\xi=
\begin{cases} 0 &\text{ if }f|_{V_\xi}=0,\\
\exp(-v_\xi(f|_{V_\xi})) &\text{ otherwise.}
\end{cases}$$
The topology on $\frX$ is such that the map $\ker$ is continuous; in fact it can be defined as the weakest topology such that:
\begin{enumerate}
	\item The kernel map is continuous.
	\item For every open affine $U\subset X$ and every $f\in \cO_X(U)$, the function
\begin{align*}
\ker^{-1}(U) &\rightarrow \bbR\\
\xi&\mapsto \norm{f}_\xi
\end{align*}
is continuous.
\end{enumerate}

There is a second map $\red:\frX\rightarrow X$, called the `reduction map' in the literature, which sends each $\xi=(V_\xi,v_\xi)$ to $\cent(v_\xi)$, the center of the valuation $v_\xi$, defined as the prime ideal in $\cO_X(U)$ formed by those $f\in \cO_X(U)$ with $\norm{f}_\xi<1$ (for any affine neighborhood $U$ of $\ker \xi$). 
The (schematic) point $\red(\xi)$ belongs to $V_\xi$, and there, it can be described as the prime ideal in $\cO_{V_\xi}(U\cap V_\xi)$ formed by those $f$ with $v_\xi(f)>0$.
The reduction map is anticontinuous \cite{Tem15}.

\begin{notation}
For every irreducible and reduced $V\subseteq X$, the notation $\ordnorm_V=(X,\ord_V)\in \frX$ will stand for the norm induced by the valuation \emph{order along} $V$, i.e., $\ker(\ordnorm_V)=0$ and for every $f\in K(X)$, 
$\norm{f}_{\ordnorm_V}=\exp(-\ord_{V}(f))$.
 	
For every $\xi=(V_\xi,v_\xi) \in \frX$, and $\ell\in[0,+\infty)$ we use the notation $\xi^\ell$ for the seminorm $\xi^\ell=(V_\xi,\ell v_\xi)$; in other words, 
\begin{gather*}
\norm{f}_{\xi^\ell}=0 \Longleftrightarrow \norm{f}_\xi=0,\\
\norm{f}_{\xi^\ell}=\norm{f}_\xi^\ell\quad \text{ if } \norm{f}_\xi \ne 0.
\end{gather*}
In particular, for $\ell=0$ we get the trivial seminorm on $V_\xi = \overline{\ker(\xi)}$, satisfying 
$$\norm{f}_{\xi^0}=\begin{cases}
0 & \text{if }f|_{V_\xi}=0,\\ 1 &\text{otherwise.}
\end{cases} $$
Additionally, by $\xi^\infty$ we mean the trivial seminorm on $\red(\xi)$, i.e. the seminorm which satisfies
$$\norm{f}_{\xi^\infty}=\lim_{\ell\to\infty}\norm{f}_\xi^\ell=
\begin{cases}
0& \text{if } f|_{\overline{\red(\xi)}}=0,\\
1& \text{otherwise,}
\end{cases}$$ 
for every $f$ regular in a neighborhood of $\red(\xi)$.
Thus, the trivial seminorm on a given subvariety $V$ shall be denoted by $\ordnorm_V^\infty$.
Note that for the trivial norm we have $\ordnorm_X=\ordnorm_X^\infty$.
\end{notation}

Given a real valuation $v$ on $K(X)$ with center $x$ on $X$, one has a real filtration by ideal sheaves $(\fra_{v,m} \subset \cO_X)_{m\in\bbR_{\ge 0}}$ defined by
\[ 
\fra_{v,m}(U)=\begin{cases}
\{f \in \cO_X(U) , v(f)\ge m\} & \text{if }x\in U,\\
\cO_X(U) & \text{ otherwise}.
\end{cases} 
\]
This filtration is multiplicative, i.e., $\fra_{v,m}\cdot\fra_{v,n}\subset\fra_{v,m+n}$ for every $m, n \in \bbR_{\ge0}$.
This generalizes immediately to the more general setting of a seminorm $\xi \in \frX$;
recall that $\xi$ can be described as a pair (kernel, valuation), $\xi=(V_\xi,v_\xi)$, with $v_{\xi}(f|_{V_\xi}) = -\log \norm{f}_\xi$. 
Then the filtration is defined as
\[ 
\fra_{\xi,m}(U)=\begin{cases}
\{f \in \cO_X(U) , \log \norm{f}_\xi\le-m\} & \text{if }\red(\xi)\in U,\\
\cO_X(U) & \text{ otherwise}.
\end{cases} 
\]
Sometimes we need to emphasize the model $X$ where these ideal sheaves live, and we shall use the notation $\fra_{\xi,m}^X$.

\subsection{Riemann-Zariski space and Shokurov \textbf{b}-divisors}
\label{sec:intro-b-divisors}

The Riemann--Zariski space of $X$ was introduced by Zariski \cite[Chapter VI,\S17]{ZS75II} as the set of all equivalence classes of valuations on $K(X)$ with a suitable structure of locally ringed space, which makes it isomorphic to the projective limit 
$$
\cX:=\varprojlim_{\pi} X_{\pi}
$$
(in the category of locally ringed spaces)
where $\pi\colon X_{\pi} \to X$ is a proper birational morphism.
Such a limit is not a scheme if $\dim X>1$.
For every point $\eta\in\cX$, the local ring $\cO_{\cX,\eta}$ is a valuation ring of $K(X)$, namely the common valuation ring of all valuations in the equivalence class $\eta$, and every valuation ring arises this way.

The birational study of divisors on models of $X$ is most conveniently tackled by using Shokurov's language of \textbf{b}-divisors on the Riemann-Zariski space, introduced in \cite{Sho03}. Let us recall some well-known facts following \cite{BdFF12}. 

The group of Weil and Cartier \textbf{b}-divisors over $X$ are respectively defined as
\begin{eqnarray*}
\mathrm{Div}(\cX)&:=& \varprojlim_{\pi} \mathrm{Div}(X_{\pi}) \quad \textrm{and}\\
\mathrm{CDiv}(\cX)&:=& \varinjlim_{\pi} \mathrm{CDiv}(X_{\pi}),
\end{eqnarray*}
where the first limit is taken with respect to the pushforward and the second with respect to the pullback.
As is usual, $\bbR$-Weil and $\bbR$-Cartier \textbf{b}-divisors are obtained by extension of scalars, tensoring these abelian groups with $\otimes_{\bbZ}\bbR$.

\begin{itemize}
    \item Given a Weil \textbf{b}-divisor $D$, the \emph{trace} (or incarnation) of $D$ on any given model $X_{\pi}$ will be denoted by $D_{\pi}:=D|_{X_{\pi}}$.
    \item We will say that a Cartier \textbf{b}-divisor $C$ is \emph{determined} on $X_{\pi}$ for a model $X_{\pi}$ over $X$ if the equality $C_{\pi'}=\phi^*C_{\pi}$ holds for every higher model $\phi\colon X_{\pi'} \to X_{\pi}$.
    In that case we also say that $C_\pi$ is a \emph{determination} of $C$.
    \item Let $X_{\pi} \overset{\pi}{\rightarrow} X$ be some model and $D_{\pi}$ a Cartier divisor on $X_{\pi}$. We denote by $\overline{D_{\pi}}$ the Cartier closure of $D_{\pi}$: the Cartier \textbf{b}-divisor obtained via pullback on higher models.
    \item Given a Weil \textbf{b}-divisor $D$, the fractional ideal sheaf on $X$ whose sections on an open set $U \subset X$ are the rational functions $f$ such that $\overline{\div(f)} + W \ge 0$ is denoted by $\cO_X(D)$.
    This sheaf is not coherent in general, but it is coherent for Cartier \textbf{b}-divisors and also for some important classes of Weil divisors, as will be seen below.
\end{itemize}

We will need the notion of numerical equivalence for \textbf{b}-divisors from \cite{BdFF12}. There are two possible natural definitions:
\begin{itemize}
\item the 1-codimensional numerical classes over $X$ 
$$N^1(\cX) := \varinjlim_{\pi}N^1(X_{\pi})$$
with respect to pulling-back. 
\item the $(n - 1)$-dimensional numerical classes over $X$ 
$$N_{n-1}(\cX) := \varprojlim_{\pi} N^1(X_{\pi})$$
with respect to push-forward.
\end{itemize}

There is a natural continuous injection $N^1(\cX)\hookrightarrow N_{n-1}(\cX)$, as well as surjections $\CDiv_{\bbR}(\cX)\rightarrow N^1(\cX)$, $\Div_{\bbR}(\cX)\rightarrow N_{n-1}(\cX)$ compatible with it, which we write simply as $D\mapsto [D]$.
The latter map is not continuous in general, but one has the following, which is often sufficient in applications:

\begin{lemma}{\cite[Lemma 1.12]{BdFF12}}\label{lem:convergence}
Let $W_j$ be a sequence (or net) of $\mathbb{R}$-Weil \textbf{b}-divisors which converges to an $\mathbb{R}$-Weil \textbf{b}-divisor $W$ coefficient-wise. If there exists a fixed finite dimensional vector space $V$ of $\mathbb{R}$-Weil divisors on $X$ such that $W_{j,X} \in V$ for all $j$ then $[W_j] \rightarrow [W]$ converges in $N_{n-1}(\cX)$.
\end{lemma}

With numerical equivalence at hand, one can deal with the notion of \emph{nefness} of Cartier and even Weil divisors:

\begin{itemize}
    \item A Cartier \textbf{b}-divisor is said to be $X$-nef or nef over $X$ if for every model $X_{\pi} \overset{\pi}{\rightarrow} X$ where it is determined and every curve $C \subseteq X_{\pi}$ that gets contracted to a point on $X$, we have $D_{\pi}\cdot C\geq 0$.
    \item A Weil \textbf{b}-divisor $D$ is said to be $X$-nef or nef over $X$ if its numerical class is the limit of a sequence or net of $X$-nef Cartier \textbf{b}-divisors. 
    Equivalently, for every smooth model $X_{\pi} \overset{\pi}{\rightarrow} X$, $D_\pi$ belongs to the closed movable cone \cite[Lemma 2.10]{BdFF12}.
    \item A Cartier \textbf{b}-divisor $D$ which is determined on the model $X_\pi$ is said to be \emph{relatively globally generated} over $X$ if so is $D_\pi$, i.e., if $\pi^*\pi_*\cO_{X_\pi}(D_\pi) \rightarrow \cO_{X_\pi}(D_\pi)$ is surjective.
\end{itemize}

\begin{definition}
Let $\fra$ be a coherent fractional ideal sheaf on $X$. If $X_{\pi} \to X$ is the normalized blow-up of $X$ along $\mathfrak{a}$, the Cartier \textbf{b}-divisor determined by $\mathfrak{a}\cdot \cO_{X_{\pi}}$ is called $Z(\mathfrak{a})$.

Let $\fra_{\bullet}$ be a graded sequence of (coherent) fractional ideal sheaves with linearly bounded denominators (see \cite[Section 2.1]{BdFF12}).
Then the limit $\lim_{m\to \infty} \frac{1}{m} Z(\mathfrak{a}_m)$ is an $\bbR$-Weil \textbf{b}-divisor called $Z(\fra_\bullet)$. 
\end{definition}

A Cartier \textbf{b}-divisor $C \in \mathrm{CDiv}(\cX)$ is of the form $Z(\mathfrak{a})$ for some fractional ideal $\fra$ if and only if $C$ is relatively globally generated over $X$ \cite[Lemma 1.8]{BdFF12}.

\begin{definition}[\cite{BdFF12}, Definition 2.3]\label{def:enf-envelope} Let $D$ be an $\mathbb{R}$-Weil divisor on some model $X_\pi$. The \emph{nef envelope} $\mathrm{Env}_\pi(D)$ is the $\mathbb{R}$-Weil \textbf{b}-divisor associated with the graded sequence $\mathfrak{a}_{\bullet}=\{\pi_ *\cO_{X_\pi}(mD)\}_{m \geq 1}$, namely $$Z(\mathfrak{a}_{\bullet}) = \lim_{m\to \infty} \frac{1}{m} Z(\mathfrak{a}_m).$$
\end{definition}

\begin{lemma}[\cite{BdFF12}, Lemma 2.11 and Corollary 2.13]\label{lem:nef-envelope-ZR} We have that 
\begin{enumerate}
\item $Z(\mathfrak{a}_{\bullet})$ is $X$-nef for every graded sequence of fractional ideals with linearly bounded denominators, and
\item $\mathrm{Env}_\pi(D)$ is the largest $X$-nef $\mathbb{R}$-Weil \textbf{b}-divisor $W$ such that $W_{\pi}\leq D$.
\end{enumerate}
\end{lemma}

Nef envelopes of \emph{arbitrary} \textbf{b}-divisors $W$ are not defined, but only because the existence of $X$-nef \textbf{b}-divisors bounded by $W$ is not guaranteed. Under this hypothesis they can be defined:

\begin{definition}[\cite{BdFF12}, Proposition 2.15 and Definition 2.16]\label{def:nef-envelope-ZR}
Let $W$ be an $\bbR$-Weil \textbf{b}-divisor. If the set of $X$-nef $\bbR$-Weil \textbf{b}-divisors $Z$ such that $Z\le W$ is nonempty, then it admits a largest element $\Env_{\cX}(W)$, called the nef envelope of $W$.
\end{definition}

\subsection{Evaluation of seminorms on Cartier \textbf{b}-divisors}

Given a seminorm $\norm{\cdot}_\xi$, let $U_\xi$ be an affine open neighborhood of $\red(\xi)$ such that $\norm{\cdot}_\xi$ is defined on $\Gamma(\cO_X,U)$. If $D$ is an effective Cartier divisor on $X$, and $f\in\Gamma(\cO_X,U)$ is an equation of $D$ in $U$, set $\norm{D}_\xi=\norm{f}_\xi$. 
It is easy to see that such a definition does not depend on the choices,
and it extends to arbitrary Cartier divisors by multiplicativity.

Since pullback by birational maps preserves \emph{values} of any $f\in\Gamma(\cO_X,U)$, it is clear that \emph{if $\xi\in \frX$ is a norm}, i.e., it has an associated valuation $v_\xi$ on $K(X)$ ($\ker \xi$ is the generic point of $X$), then $\norm{\cdot}_\xi$ can be evaluated on all Cartier \textbf{b}-divisors, by evaluating on an arbitrary representative. 

For seminorms which are not necessarily norms, some care is needed in evaluating them on higher models:

\begin{lemma}
 Let $D$ be a Cartier \textbf{b}-divisor, determined on some model $X_\pi\overset{\pi}{\rightarrow} X$, and let $U_\pi\subset X$ be an open set over which $\pi$ is an isomorphism.

 Consider the open set $\frU_\pi\subset\frX$ of all seminorms with $\ker(\xi)\in U_\pi$.
 For every $\xi\in \frU_\pi$, 
 $\norm{D}_\xi$ is well defined, and the map $\frU_\pi \to \bbR$ defined by
 $\xi \mapsto \norm{D}_\xi$
 is continuous. 
\end{lemma}
\begin{proof}
We assume without loss of generality that $D$ is effective.

Berkovich analytification is functorial, so $\pi$ determines an analytic morphism $\pi^\beth:\frX_\pi \rightarrow \frX$ from the analytification of $X_\pi$ to $\frX$, which restricts to an isomorphism over $\frU_\pi$. 
Therefore all seminorms in $\frU_\pi$ can be evaluated on Cartier divisors of $X_\pi$, such as $D_\pi$, i.e., $\norm{D}_\xi = \norm{D_\pi}_{(\pi^{\beth})^{-1}(\xi)}$ is well defined for every $\xi\in \frU_\pi$.

Now fix $\xi_0\in\frU_\pi$, and let $P$ be the center of $\xi_0$ on $X_\pi$, i.e., $P=\red_{X_\pi}(\xi_0)$. Since the map $\red_{X_\pi}$ is anticontinuous, the set $\frU_0$ of all $\xi \in \frU_\pi$ such that $\red_{X_\pi}(\xi)\in \overline{\red_{X_\pi}(\xi_0)}$ is an open neighborhood of $\xi_0$ in $\frU_\pi\subset\frX$, and we need only show continuity of $\norm{D}_\xi$ for $\xi \in \frU_0$.
Let $U_1,\dots, U_k$ be open affine neighborhoods of $\red_{X_\pi}(\xi_0)$ covering $\overline{\red_{X_\pi}(\xi_0)}$, and $f_i\in \Gamma(\cO_{X_\pi},U_i)$ equations of $D$ in $U_i$. 
In the (closed) subset $\mathfrak{Z}_i\subset \frU_0$ consisting of the seminorms with center in $U_i$, the function $\xi\mapsto\norm{D}_\xi=\norm{f_i}_\xi$ is continuous by the definition of the topology on $\frX$.
On the other hand, on overlaps $U_i\cap U_j$ we have $\norm{f_i}_\xi=\norm{D}_\xi= \norm{f_j}_\xi$ because the norm of $D$ is well defined (explicitly, the quotient $f_i/f_j$ is invertible, so it does not belong to $\red_{X_\pi}(\xi)$ for any $\xi \in \mathfrak{Z}_i\cap\mathfrak{Z}_j$, i.e., its seminorm is 1).
We conclude that the functions $f_i$ patch to give a continuous function $\xi\mapsto\norm{D}_\xi$ on $\frU_0$.
\end{proof}

\begin{notation}
For every regular function $f\in \Gamma(\cO_{X},U)$ on some open subset $U\subset X$, write $\div_+(f)$ for the effective divisor obtained as the closure of $V(f)\subset U$. 
Then for every seminorm $\xi$ such that $f$ is regular at $\red \xi$, we have $\norm{\div_+(f)}_\xi=\norm{f}_\xi\le 1$.
If $f$ is not regular at $\red \xi$, then by definition $\norm{\div_+(f)}_\xi=\norm{f'}_\xi$ where $f'$ is an equation for $\div_+(f)$ regular at $\red \xi$, and so in all cases $\norm{\div_+(f)}_\xi\le 1$ (equivalently, $v_\xi(\div_+(f))\ge 0$).
\end{notation}

\subsection{Quasi-monomial valuations and seminorms}
\label{sec:quasi-monomial}
Fix a model $X_\pi$ of $X$ and a point $p\in X_\pi$ such that the local ring $\cO_{X_\pi,p}$ is regular of dimension $c\le \dim X$. 
Let $z_j\in \cO_{X_\pi,p}$, $j=1,\dots, c$ be a regular system of parameters, and let $D_j=\div_+(z_j)$. By Cohen's theorem, $\widehat{\cO}_{X_\pi,p}\cong \kappa(p)[[z_j]]$, where $\kappa(p)$ is the residue field of $\cO_{X_\pi,p}$, so every $f\in \cO_{X_\pi,p}$ has a unique expansion as $f=\sum_{\alpha \in \bbN^c} a_\alpha z^{\alpha}$.
To every weight $w=(w_j)\in \bbR_{\ge0}^c$ one associates the monomial valuation $v_w$ of $\cO_{X_\pi,p}$ determined by
\[
v_w\left(\sum_{\alpha \in \bbN^c} a_\alpha z^{\alpha}\right) = 
\min\left\{ w\cdot \alpha \,|\, a_\alpha\ne 0 \right\},
\]
which extends to a valuation of $K(X)$, and determines a norm $\norm{f}_w=\exp (-v_w(f))$.
We may further allow $w\in (\bbR_{\ge0} \cup \{+\infty\})^c$, in which case the prime ideal $\fra_{w,\infty}=(z_j|w_j=+\infty)$ consists of all $f\in {\cO}_{X_\pi,p}$ with $v_w(f)=+\infty$, $v_w$ determines a valuation on the quotient ${\cO}_{X_\pi,p}/\fra_{w,\infty}$, and $\norm{f}_w=\exp (-v_w(f))$ is a seminorm.
In all cases, the center of $v_w$ is the prime ideal $(z_j)_{w_j>0}$.
The map $(\bbR_{\ge0} \cup \{+\infty\})^c \rightarrow \frX$ determined by $w \mapsto \norm{\cdot}_w$ is a homeomorphism onto its image $\Delta_{\pi,p,z}$, and one has $\frX=\overline{\bigcup_{\pi,p,z}\Delta_{\pi,p,z}}$ (\cite{Ber90,Thu07,JM12}).

In fact, one can more precisely describe $\frX$ as an inverse limit of such subsets $\Delta_{\pi,p,z}$.
The map
\begin{equation}\label{eq:retraction}
\begin{aligned}
 ev_{z}:\frX &\rightarrow(\bbR_{\ge 0}\cup\{+\infty\})^c\\
 \xi &\mapsto (-\log \norm{D_1}_\xi,\dots,-\log \norm{D_c}_\xi)
\end{aligned}\end{equation}
is a continuous retraction and satisfies $\norm{f}_{ev_{z}(\xi)}\ge \norm{f}_\xi$ for every $f\in \mathcal{O}_{X_\pi,p}$. 
For every log-smooth pair $(X_\pi,D)$ over $X$, 
call the pair $(p,z)$ $D$-adapted if each $D_j=\div_+(z_j)$ is one of the irreducible components of $D$.
Then the topology of the union
$$
\Delta_{\pi,D}=\bigcup_{(p,z)\ D\text{-adapted}} \Delta_{\pi,p,z} \subset \frX
$$
is the topology of a simplicial cone complex obtained by gluing the cones $\Delta_{\pi,p,z}\simeq (\bbR_{\ge0} \cup \{+\infty\})^c$, and $\frX$ is the inverse limit of all such simplicial complexes via the retraction map. This is explained in all detail for the case $w\in \bbR_{\ge0}^c$ of valuations in \cite{JM12}; the generalization to $w\in (\bbR_{\ge0} \cup \{+\infty\})^c$ and seminorms is straightforward.

Valuations $v_w$ and seminorms $\norm{\cdot}_w$ as above which are monomial \emph{on some model} $X_\pi$ are generally called \emph{quasi-monomial}.
It is known by \cite[Proposition 2.8]{ELS03} that quasi-monomial valuations can be algebraically characterized as those valuations for which Abhyankar's inequality
$$
\trdeg v + \ratrk v \le \dim X
$$
is in fact an equality (where $\trdeg$ stands for the transcendence degree of the residue field of $v$, and $\ratrk$ stands for the rational rank of the value group of $v$).
For this reason, quasi-monomial valuations are also called Abhyankar valuations in the literature.

\subsection{Noether lemma}
When dealing with \textbf{b}-divisors, the role of \emph{primes} is most naturally played by those associated to divisorial valuations.
 The center of a divisorial valuation on every sufficiently high model is a divisor; when passing to another model, the centers are related by strict transforms and push-forwards, and thus give rise to a Weil \textbf{b}-divisor:

 \begin{notation} Let $X_{\pi}$ be any birational model of $X$. Let $D$ be a prime divisor on $X_{\pi}$. We will denote by $\widetilde{D}$ the \textbf{b}-divisor obtained from $D$ via pushforward and strict-transform.
 \end{notation}

 A Weil \textbf{b}-divisor $\widetilde{D}$ is called \emph{prime} if it is the \textbf{b}-divisor  determined by strict transforms and pushforwards of a prime divisor $D$ on some model $\pi:X_\pi\rightarrow X$.
 A prime \textbf{b}-divisor is not Cartier, but it can be easily seen as a limit of Cartier \textbf{b}-divisors.

\begin{proposition}\label{pro:separate}
	Let $v$ be a rank 1 valuation of $K(X)$, and let $\widetilde{D}$ be a prime 
    Weil \textbf{b}-divisor, determined by a prime divisor $D$ on the model $X_\pi$. 
	If $v$ is not equivalent to the divisorial valuation $w$ associated to $D$, there exists a model $\varpi:X_\varpi\rightarrow X_\pi$ such that the center of $v$ in $X_\varpi$ is not contained in $\widetilde{D}_\varpi$.
\end{proposition}
\begin{proof}
	The valuation rings of $K(X)$ strictly contained in $\cO_w$ are valuation rings of valuations composed with $w$, in particular they have rank greater than one. 
	So, if $v$ is not equivalent to $w$, $\cO_v$ is not contained in $\cO_w$. 
	Let $f\in \cO_v\setminus \cO_w$.
	Then 
	$$\mathcal{U}_f=\{\mu \text{ valuation of }K(X)\,|\,\mu(f)\ge 0\}$$
	is open in $\cX$, it contains $v$ and it does not contain $w$. 
	Since $\cX$ is the limit of all models $X_\varpi$, there is one such model where the center of $v$ has a neighborhood that does not contain the center of $w$, which is the generic point of $\widetilde{D}_\varpi$.
\end{proof}

Suppose we are given a sequence of blowups centered at smooth irreducible subvarieties,
\[X_n \overset{bl_{V_n}}{\longrightarrow}X_{n-1} \overset{bl_{V_{n-1}}}{\longrightarrow}
\dots  \overset{bl_{V_2}}{\longrightarrow}X_1 \overset{bl_{V_1}}{\longrightarrow} X_0=X .\]
For $i=1,\dots, n$, denote $E_i$ the exceptional divisor of the blowup $bl_{V_i}$ above the subvariety $V_i\subset X_{i-1}$, and for every effective divisor $D$ on $X$, denote $\tilde D_i$ the strict transform of $D$ on $X_i$, and $m_i(D)$ the multiplicity of $\tilde D_i$ along $V_i$.

\begin{lemma}[{Noether equality, \cite[II.2]{Nun04}}]\label{lem:noether}
	Let $D$ be an effective divisor on $X$, $v$ a valuation of $K(X)$.
	For every sequence of blowups with smooth centers as above,
	\[v(D)=\sum_{i=1}^{n} m_i(D)v(E_i)+ v(\tilde D_n).\]
	If $v$ is a real valuation, then there is a sequence of blowups with smooth centers such that $v(\tilde D_n)=0$.
\end{lemma}

\begin{proof}
	Write the pullback of $D$ to $X_n$ as a combination of the pullbacks of the $E_i$ and $\tilde D_n$.
	For the second claim, use the previous proposition.
\end{proof}

\subsection{Multiplier ideals and Arnold multiplicity}

The asymptotic behavior of valuation ideal sheaves is conveniently controlled by so-called asymptotic multiplier ideals. 
We recall their definition and main properties from \cite[Part Three]{LazP2}. 

Let $\fra\subseteq\cO_X$ be a non-zero ideal sheaf, and $c>0$ a rational number. 
Fix a log-resolution $\pi:X' \rightarrow X$ of $\fra$ with $\fra\cdot\cO_{X'}=\cO_{X'}(Z(\fra))$.
The \emph{multiplier ideal} $\cJ(\fra^c)$ associated to $c$ and $\fra$ is defined as
$$
\cJ(\fra^c)=\pi_*\cO_{X'}(K_{X'/X}+\lceil c\cdot Z(\fra)\rceil)
$$
where $K_{X'/X}$ stands for the relative canonical divisor and $\lceil c\cdot Z(\fra) \rceil$ is the least integral divisor $D\ge c\cdot Z(\fra)$\footnote{We depart from the notational convention used in \cite[Definition 9.2.3]{LazP2} which uses the effective divisor $F=-Z(\fra)$ and round-down; our choice looks for consistency with the usual convention when working with Nef envelopes: $Z(\fra)$ is $X$-nef whereas $F$ is anti-$X$-nef.}. 
The resulting ideal is independent of the choice of log-resolution $\pi$, it contains $\fra$, and the subadditivity relation 
\begin{equation}\label{eq:subadditivity}
 \cJ(\fra^{\ell c})\subseteq \cJ(\fra^c)^\ell   
\end{equation}
holds for every $\ell \in \bbN$ \cite[Theorem 9.5.20]{LazP2}.

The \emph{Arnold multiplicity} of an ideal $\fra$ at $p\in V(\fra)$ can be defined as
$$
\operatorname{Arn}_p(\fra) = \left(\inf\{c>0 \mid \cJ(\fra^c)\subset\frm_p\}\right)^{-1},
$$
i.e., it is the inverse of the \emph{log-canonical threshold} (in addition to \cite{LazP2},
see \cite{Kol97}, \cite{JM12}).

For $\fra_\bullet$ a graded sequence of (coherent) ideal sheaves, such as the valuation ideals we are mostly concerned with, one has $\cJ(\fra_m)\subseteq \cJ(\fra_{pm}^{1/p})$ for every natural $p$ \cite[Lemma 11.1.4]{LazP2}.
By noetherianity then, there is a unique maximal element in the set of ideals $\{\cJ(\fra_{pm}^{1/p})\}_{p\in\bbN}$, which is called the $m$-th asymptotic multiplier ideal $\frj_m=\frj_m(\fra_\bullet)$. 
Its essential property is that for every $m$ and every natural number $\ell$, 
\begin{equation}
    \label{eq:asymptotic_ideal_inclusion}
    \fra_m^\ell \subseteq \fra_{\ell m} \subseteq \frj_m^\ell,
\end{equation}
which is a consequence of subadditivity \cite[Theorem 11.2.3]{LazP2}: $\frj_{s+t}\subseteq\frj_s\cdot\frj_t$.
The Arnold multiplicity of the graded sequence is defined as $$\operatorname{Arn}_p(\fra_\bullet):=
\inf_{m \in\bbN}\frac{1}{m}\operatorname{Arn}_p(\fra_m)=
\lim_{m\to\infty} \frac{1}{m}\operatorname{Arn}_p(\fra_m).$$

The divisor $K_{X'/X}+\lceil c\cdot Z(\fra)\rceil$ is not necessarily $X$-nef,
hence we will be interested in the ``multiplier ideal divisor''
$$
Z(\fra^c):=Z(\cJ(\fra^c)) = \Env_\pi(K_{X'/X}+\lceil c\cdot Z(\fra)\rceil)
$$
seen as a Cartier \textbf{b}-divisor. 
Again, for our purposes, we are mainly interested in the asymptotic counterparts; 
given a graded sequence of ideals $\fra_\bullet$, we define the $m$-th asymptotic 
multiplier ideal divisor as the Cartier \textbf{b}-divisor
$$
\mathring Z_m(\fra_\bullet):=Z(\frj_m(\fra_\bullet)). 
$$
By \eqref{eq:asymptotic_ideal_inclusion}, these \textbf{b}-divisors satisfy
\begin{equation}
\label{eq:above-below}
Z(\fra_m) \le \frac{Z(\fra_{\ell m})}{\ell}\le \mathring Z_m(\fra_\bullet) 
\end{equation}
for every $m>0$ and every natural $\ell$.

\medskip 

Recall that a \emph{complete} ideal is an
ideal determined by the values $\mathbf{m}=(m_i)_{i\in I}$ of its elements at some set of valuations $v_i$
$$ 
\fra_{\mathbf{m}} = \{f \in \cO_X \mid v_i(f)\ge m_i \,\forall i \in I\}.
$$
This is equivalent to the ideal being integrally closed \cite[Chapter 6]{HS06}.
\begin{proposition}\label{pro:envelope-under-birational}
If $\phi:X\rightarrow Y$ is a birational morphism, and $\fra=\fra_{\mathbf{m}}$ is
a complete ideal, then $Z(\phi_*(\fra))=\Env_{\cY}(Z(\fra))$.
\end{proposition}
Note that, since $Z(\fra)$ is a Cartier \textbf{b}-divisor determined on the normalized blowup 
$\pi:X_\pi\rightarrow X$ of $X$ along $\fra$, one has 
$\Env_{\cY}(Z(\fra))=\Env_{\phi\circ\pi}((Z(\fra))_\pi)$.
\begin{proof}
Note first that $\phi_*(\fra)$ is also a valuation ideal
(given by the same set of valuations on $K(Y)=K(X)$),
and $\phi^* \phi_*(\fra) \subseteq \fra$.

    Let $\pi_Y:Y_{\mathbf{m}}\rightarrow Y$ be the normalized blowup of $Y$ along $\phi_*(\fra)$, and let $\pi':X_{\mathbf{m}}\rightarrow X_\pi$ be a birational model such that the induced map $\phi_{\mathbf{m}}:X_{\mathbf{m}}\rightarrow Y_{\mathbf{m}}$ is a morphism.
    For simplicity write also $\pi_X=\pi\circ\pi'$.
    
    By the definitions then, the Cartier \textbf{b}-divisors $Z(\phi_*(\fra))$ and $Z(\fra)$ are defined on $Y_{\mathbf{m}}$ and $X_{\mathbf{m}}$ respectively, and
	\begin{align*}
	\fra&={\pi_X}_*(\cO_{X_{\mathbf{m}}}(Z(\fra))_{X_{\mathbf{m}}}),\\
	\phi_*(\fra)&={\pi_Y}_*(\cO_{Y_{\mathbf{m}}}(Z(\phi_*(\fra)))_{Y_{\mathbf{m}}}).
	\end{align*}
    We will prove that $\phi_{\mathbf{m}}^*(Z(\phi_*(\fra))_{Y_{\mathbf{m}}}) = \Env_{\pi_Y\circ\phi_{\mathbf{m}}}(Z(\fra)_{X_{\mathbf{m}}})$.
    This will show that the two \textbf{b}-divisors in the claim are equal, as they are defined on $X_\mathbf{m}$ and agree there, by commutativity of the diagram
$$	\begin{tikzcd}
X_{\mathbf{m}}\arrow[r,"\phi_\mathbf{m}"] \arrow[d,swap,"\pi'"]& Y_{\mathbf{m}} 
\arrow[dd,swap,"\pi_Y"] \\
X_{\pi}  \arrow[d,swap,"\pi"]&\\
X\arrow[r,"\phi"]&Y 
\end{tikzcd}$$

	Since $\phi^* \phi_*(\fra) \subseteq \fra$, we deduce that $\phi_{\mathbf{m}}^*(Z(\phi_*(\fra))_{Y_{\mathbf{m}}}) \leq Z(\fra)_{X_{\mathbf{m}}}$. Since $Z(\phi_*(\fra))_{Y_{\mathbf{m}}}$ is $Y$-nef, it follows by \cite[Corollary 2.13]{BdFF12} that $\phi_{\mathbf{m}}^*(Z(\phi_*(\fra))_{Y_{\mathbf{m}}}) \leq \Env_{\pi_Y\circ\phi_{\mathbf{m}}}(Z(\fra)_{X_{\mathbf{m}}})$.
	
	To see the converse inequality $\phi_{\mathbf{m}}^*(Z(\phi_*(\fra))_{Y_{\mathbf{m}}}) \geq \Env_{\pi_Y\circ\phi_{\mathbf{m}}}(Z(\fra)_{X_{\mathbf{m}}})$, it will be enough to see that, for every Cartier divisor $D$ on $Y$ such that $(\phi\circ\pi_X)^*(D)+\Env_{\pi_Y\circ\phi_{\mathbf{m}}}(Z(\fra)_{X_{\mathbf{m}}})\ge 0$, the inequality $(\phi\circ\pi_X)^*(D)+\phi_{\mathbf{m}}^*(Z(\phi_*(\fra))_{Y_{\mathbf{m}}}) \geq 0$ is also satisfied.
	
	Now for every $f\in \cO_Y$ such that $(\phi\circ\pi_X)^*(\div(f))+ \Env_{\pi_Y\circ\phi_{\mathbf{m}}}(Z(\fra)_{X_{\mathbf{m}}})\geq 0$ we get $(\phi\circ\pi_X)^*(\div(f))+ Z(\fra)_{X_{\mathbf{m}}}\geq 0$, so 
 $f\in \phi_*\fra$ which implies $v_i(\phi^*f)\ge m_i$ for all $i\in I$ and therefore
 $v_i(f)\ge m_i$. Thus $(\pi_Y)^*(\div(f))+Z(\phi_*(\fra))_{Y_{\mathbf{m}}}\ge 0$ and $(\phi\circ\pi_X)^*(\div(f))+\phi_{\mathbf{m}}^*(Z(\phi_*(\fra))_{Y_{\mathbf{m}}})\ge 0$ as claimed.
\end{proof}

\begin{corollary}\label{cor:Zxi-under-birational}
If $X\rightarrow Y$ is a birational morphism between smooth varieties, 
and $\fra_\bullet$ is a graded sequence of valuation ideal sheaves on $X$, 
then $Z(\phi_*(\fra_m))=\Env_\cY(Z(\fra_m))$ and
$\mathring Z_m(\phi_*(\fra_\bullet))\ge\Env_\cY(\mathring Z_m(\fra_\bullet))$ for every $m$.
\end{corollary}
\begin{proof}
    The equality $Z(\phi_*(\fra_m))=\Env_\cY(Z(\fra_m))$ is immediate from the proposition,
    as $\fra_m$ is a complete ideal for every $m$.
    For multiplier ideal divisors, we have
    \begin{gather*}       
    \mathring Z_m(\phi_*(\fra_\bullet))=Z(\frj_m(\phi_*(\fra_\bullet)))=
    \Env_{\phi\circ\pi}(K_{X_m/Y}+\lceil (1/p) \cdot Z(\phi_*(\fra_{pm}))_{X_m}\rceil) = \\
    \Env_{\phi\circ\pi}(Z(\frj_m(\fra_\bullet))+\pi^*(K_{X/Y})).
    \end{gather*}
    Since $K_{X/Y}\ge 0$, it follows that 
    $$
    \mathring Z_m(\phi_*(\fra_\bullet))\ge \Env_{\phi\circ\pi}(Z(\frj_m(\fra_\bullet)))=
    \Env_\cY(\mathring Z_m(\fra_\bullet))
    $$
    as claimed.
\end{proof}

\section{The \textbf{b}-divisor associated to a seminorm}
\label{sec:d-xi}

In this section we define the map  $D:\frX\to \Div\cX$ which is central in our work and we study its formal properties. In later sections it will be used to connect valuative properties of divisors with positivity properties.

\subsection{\texorpdfstring{$X$}{X}-nef \textbf{b}-divisor associated to an analytic point \texorpdfstring{$\xi\in\frX$}{ξ}}

Recall from section \ref{sec:intro-Berkovich} that we associate to every seminorm $\xi \in \frX$ a real filtration by ideal sheaves $(\fra_{\xi,m} \subset \cO_X)_{m\in\bbR}$  defined by
\[ 
\fra_{\xi,m}(U)=\begin{cases}
\{f \in \cO_X(U) , \log \norm{f}_\xi\le-m\} & \text{if }\red(\xi)\in U,\\
\cO_X(U) & \text{ otherwise}.
\end{cases} 
\]
Consider the coresponding Cartier divisors $Z(\fra_{\xi,m})$ as defined by Boucksom--de-Fernex--Favre \cite{BdFF12} (see section \ref{sec:intro-b-divisors}).

\begin{proposition}\label{pro:existence-nef-envelope}
	The limit $\lim_{m\to\infty} (1/m)Z(\fra_{\xi,m})$ exists as an $\bbR$-Weil \textbf{b}-divisor.
\end{proposition}
\begin{proof}
	This proposition is a slight generalization of  \cite[Proposition 2.1]{BdFF12}, 
    which deals with filtrations indexed by $\bbN$,  to our setting in which the 
    filtration is indexed by $m\in\bbR_{>0}$. 
    In this setting, Example 2.7 of \cite{JM12} proves that for every divisorial valuation
    $\ord_E$, the limit $\lim_{m\to\infty} (1/m)\ord_E Z(\fra_{\xi,m})$ exists;
    this is enough to show the result
    (for the topology of pointwise convergence).
\end{proof}

\begin{notation}\label{def:d-xi}
    In the sequel, we keep the notation $D_\xi$ (or $D_\xi^X$ if need be) for the Weil \textbf{b}-divisor $\lim_{m\to\infty} (1/m)Z(\fra_{\xi,m})$.
\end{notation}

\begin{remark}   
	Note that, even though valuations and filtrations by ideals indexed by ordered groups larger than $\bbR$ have been studied and are relevant in other contexts, the result does not generalize to them.
	In this sense, Proposition \ref{pro:existence-nef-envelope} is the most general version possible.
    
    To the best of our knowledge, the \textbf{b}-divisor $D_\xi$ has not been explicitly 
    introduced in the literature. 
    However, Jonsson and Musta\c{t}\u{a} do use the notation $\ord_E(\fra_{\bullet})$ for what we call $\ord_E(D_\xi)$ and so they implicitly deal with $D_\xi$, as
    the data of all these numbers is equivalent to the datum of $D_\xi$.
    Therefore several results in \cite{JM12} useful to us can be aptly expressed in the 
    language of \textbf{b}-divisors. We shall give such rephrasings in the sequel.    
\end{remark}

\begin{proposition}[{\cite[Proposition 2.13]{JM12}}]
    For every $\xi$, 
    $$
    D_\xi = \lim_{m\to \infty} \frac{\mathring{Z}_m(\fra_{\xi,\bullet})}{m}\, .
    $$
\end{proposition}

\begin{proposition}\label{pro:above_below_limit}
	For every $m\ge 0$, the following inequalities of \textbf{b}-divisors hold:
    $${Z(\fra_{\xi,m})}\le mD_\xi\le \mathring{Z}_m(\fra_{\xi,\bullet})\le0\, .$$
\end{proposition}
\begin{proof}
Immediate from \eqref{eq:above-below} taking limits for $\ell\to\infty$.
\end{proof}

\begin{remark}\label{rem:trivial_divisor}
	If $\xi=\ordnorm_X$ is the trivial norm on $X$, then the filtration $\fra_{\xi,\bullet}$ is also trivial, i.e., $\fra_{\xi,m}=0  \ \forall m$, and $D_\xi$ is not defined. We agree to write $D_{\ordnorm_X}=-\infty$, meaning that $D_{\ordnorm_X}<D$ for every \textbf{b}-divisor $D$ on $X$.
\end{remark}

\begin{remark}\label{rem:rescale-divisor}
	By definition, for every $\ell\in \bbR_{>0}$ we have $\fra_{\xi^\ell,m}=\fra_{\xi,m/\ell}$, so 
	$$D_{\xi^\ell}=\lim_{m\to\infty}\frac{1}{m}Z(\fra_{\xi,m/\ell})=\frac{1}{\ell}D_\xi.$$
\end{remark}

\begin{example}\label{exa:divisorial-envelope}
	Let $\widetilde{P}$ be a prime Weil b-divisor on $X$, determined by a prime divisor $P_\pi\subset X_\pi$. Let $v_P$, $\theta_P$ stand for the corresponding divisorial valuation and norm. 
	Then the associated divisor is a nef envelope: $D_{\theta_P}=\Env_\pi(-P_\pi)=\Env_{\cX}(-\widetilde{P})$ (see Definition \ref{def:enf-envelope}). Thus the \textbf{b}-divisor associated to a seminorm generalizes the nef envelope in the case of prime divisors. 
\end{example}

\begin{proposition}
\label{pro:coefficients-dxi}
	Let $P$ be a prime \textbf{b}-divisor on $X$, determined by an irreducible divisor $P_\pi\subset X_\pi$ on a suitable model. Then for every $\xi\in\frX$,
	\begin{enumerate}
        \item \label{part:strict_seminorm_zero}If $\xi\in \frX$ is not a norm, or equivalently, if $\xi=(V_\xi,v_\xi)$ with $V_\xi$ a proper subvariety of $X$, then $D_\xi= 0$.
	    \item \label{part:coefficients-dxi}$\ord_P(D_\xi)=-\inf_{f\in\cO_{X,\red(\xi)}}
	    \frac{\ord_{P_\pi}(\pi^*(f))}{-\log\norm{f}_\xi}$. 
	    \item If $\ord_{P_\pi}(D_{\xi,\pi})\ne 0$, then $\pi(P_\pi)\subseteq \overline{\red(\xi)}\subset X$.
	    \item $\ord_{\tilde{P}_\pi}(D_{\xi,X}) \ge -\frac{1}{v_\xi(P_\pi)}$.
	\end{enumerate}
\end{proposition}
\begin{proof}
    To prove the first claim, assume $\xi=(V_\xi,v_\xi)$ with $V_\xi$ a proper subvariety of $X$ and let $U$ be an affine open neighborhood of the center $\red(\xi)$. Note that, since the center of $\xi=(V_\xi,v_\xi)$ belongs to $V_\xi$, we have $U\cap V_\xi \ne \emptyset$. Hence there exist nonzero $f$ in $\mathscr{I}_{V}(U)\subset\cO_X(U)$, the ideal of $V$ in $U$, and such an $f$ satisfies $\norm{f}_\xi=0$. 
	Denote $D$ the effective Cartier divisor on $X$ obtained as the closure of $Z(f)$.
	Now for every $m$, $f\in \fra_{\xi,m}(U)$ and hence $\pi_m^*(D)+Z(\fra_{\xi,m})\ge 0$. Taking limits we obtain $\overline{D}+m D_\xi \ge 0$, or $0\ge D_\xi \ge -(1/m)\overline{D}$ for every $m$, and hence $D_\xi=0$.

    The second claim, in the case that $\xi$ is not a norm, follows immediately
    from the first one, as the infimum in the formula is computed as zero by those
    $f$ with $\norm{f}_\xi=0$.
    On the other hand, if $\xi$ is a norm, then because of multiplicativity of the filtration, the second claim is just a rewriting of the definition.
    
	For every natural number $m$ and any $f$ in the ideal 
	$$ \pi_*\cO_{X_\pi,P_\pi}\left(-\frac{m}{v_\xi(P_\pi)}P_\pi\right)\subset \cO_{X,\red(\xi)},$$
	$v_\xi(f)\ge m$, so this ideal is contained in $\fra_{\xi,m}$; the two last claims follow.
\end{proof}

\subsection{Characterization of nonvanishing \texorpdfstring{$D_\xi$}{valuative divisor}}

The \textbf{b}-divisor $D_\xi$ can be trivial for non-trivial seminorms, as we
saw in part \ref{part:strict_seminorm_zero} of \ref{pro:coefficients-dxi}; 
in this subsection we show further examples of the phenomenon and give a 
characterization of the underlying seminorms in therms of the Arnold multiplicity, 
inspired in the work of Jonsson-Musta\c{t}\u{a} \cite{JM12}.
In section \ref{sec:characterization-b-divisorial} we will show that in many cases 
of interest, $\xi$ can be recovered from $D_\xi$, and we
conjecture that this happens whenever $D_\xi\ne 0$.

\begin{remark}
	For every $\xi\in\frX$ distinct from the trivial norm, $\xi^\infty$ is the trivial seminorm on a strict subvariety of $X$, and so it is not a norm; hence Proposition \ref{pro:coefficients-dxi}, part \ref{part:strict_seminorm_zero} implies that the equality $D_{\xi^\ell}=\frac{1}{\ell}D_\xi$ (Remark \ref{rem:rescale-divisor}) holds for $\ell=\infty$ as well. 
	On the other hand, for $\ell=0$, recall that the power $\xi^0$ is the trivial seminorm on $\ker(\xi)$.
	If $\xi$ is not a norm, all involved divisors are trivial by the same proposition; if it is a norm then by definition $D_{\xi^0}=-\infty$, which can be seen as the limit of $\frac{1}{\ell}D_\xi$ when $\ell$ tends to zero by virtue of $D_\xi$ being anti-effective.
\end{remark}

\begin{remark}
	The converse to Proposition \ref{pro:coefficients-dxi}, part \ref{part:strict_seminorm_zero} is not true, even on surfaces.
	To see an example, let $(x,y)$ be a system of parameters of the local ring $\cO_{X,p}$ where $X$ is a smooth projective algebraic surface.
	The formal power series ring $k[[x,y]]$ can be identified with the completion of $\cO_{X,p}$. 
	Let $s(x)\in k[[x]]$ be a non-algebraic series, so that there is no polynomial $P(x,y)$ with $P(x,s(x))=0$ (since $K(X)$ is a finite algebraic extension of $k(x,y)$, there is no $f\in \cO_{X,p}$ with $f(x,s(x))=0$ either). 
	Consider the valuation on $\cO_{X,p}$ determined by $v(f)=\ord_x(f(x,s(x)))$, which extends to a valuation of $K(X)$ and determines a multiplicative norm $\norm{f}_\xi=exp(-v(f))$ and a point $\xi$ of $\frX$.
	It is not hard to check that $\ker(\xi)$ is the generic point of $X$, $\red(\xi)=p$, and $D_\xi=0$.

    Indeed, consider every truncation $f_k$ of $y-s(x)$ at order $k$; it is obvious that $\lim \norm{f_k}_\xi = 0$ whereas it is not hard to see that for every fixed $P$, $\ord_P(f_k)$ stabilizes for $k\gg 0$. Thus fer each $\varepsilon$, there exists $k\gg 0$ such that $\norm{f_k}_\xi\le \exp(-\ord_P(f_k)/\varepsilon)$. 
    Therefore, for every prime \textbf{b}-divisor $P$, $\ord_P(D_\xi)=0$ by Proposition \ref{pro:coefficients-dxi}.
\end{remark}

A complete characterization of the valuations which lead to zero or nonzero $D_\xi$
can be given using the Arnold multiplicity. 

\begin{proposition}
\label{pro:characterization-nonzero-dxi}
    If $\xi\in \frX$ is a norm, or equivalently, if $\xi=(X,v_\xi)$, then $D_\xi \ne 0$ if and only if $\operatorname{Arn}_p(\fra_{\xi,\bullet})>0$ for $p=\red_X(\xi)$.
\end{proposition}
\begin{proof}
    By definition, for each prime divisor $P$ on some model $X_\pi$ above $p$,
    $\ord_P(D_\xi)$ equals Jonsson--Musta\c{t}\u{a}'s asymptotic order
    $\ord_P(\fra_{\xi,\bullet})$. Then the result is a particular case of
    \cite[Corollary 6.10]{JM12}.
\end{proof}

\begin{corollary}
    If $v$ is a valuation of $K(X)$ such that the log-discrepancy $A_X(v)$ is finite, 
    then the associated norm $\xi=(X,v_\xi)$ satisfies $D_\xi\ne 0$. 
\end{corollary}
\begin{proof}
Let $p=\red_X(\xi)$. By \cite[Proposition 5.10]{JM12}, for every $m>0$ it holds
$$
\fra_{\xi,m}\subseteq \frm_p^{\lceil m/A(v)\rceil}
$$
so $\ord_{\tilde E_p}(D_\xi)\le -1/A(v)<0$ where $\tilde E_p$ is the prime divisor
of the exceptional divisor of blowing up $p$. 
\end{proof}

Note that this applies to Abhyankar (or quasi-monomial) valuations 
\cite[Proposition 5.1]{JM12}, so we get the following:

\begin{corollary}
    If $v$ is an Abhyankar valuation of $K(X)$ and $\xi=(X,v_\xi)$ is the associated norm, then $D_\xi\ne 0$. 
\end{corollary}

\subsection{Nonvanishing and \textbf{b}-divisorial valuations}
\label{sec:characterization-b-divisorial}

If $v$ is a valuation on $X$, and $\xi=(X,v)$ is the associated norm, the valuation of a rational function, or more generally of a Cartier divisor $D$, can often be replaced by the ``order of vanishing'' along $D_\xi$, i.e., the number
$$
\ord_{D_\xi}(D):=\sup\{t\in \bbR \,|\,\overline{D}+tD_\xi\ge 0\}.
$$
This will provide a translation of valuative properties of Cartier divisors into positivity properties, generalizing the behavior of divisorial valuations.

 \begin{proposition}\label{pro:seminorm_bound_positivity}
	Let $D$ be an effective Cartier divisor on $X$, and $\xi\in\frX$.
	Then $$\overline{D}-\log\norm{D}_\xi \cdot D_\xi\ge 0.$$
	In particular, $-\log\norm{D}_\xi \le \ord_{D_\xi}(D)$.
\end{proposition}
\begin{proof}
	The claim is empty if $D_\xi=0$, so by the previous proposition, we can assume that $\xi=(X,v_\xi)$ is given by a valuation $v_\xi$ of $K(X)$; write $m=-\log\norm{D}_\xi=v_\xi(D)\in \bbR$.
	
	For every affine $U\subset X$, and every equation $f\in\cO_X(U)$ of $D\cap U$, we have $f\in \fra_{\xi,m}(U)$, and therefore on the normalized blowup $X_{\pi_m}$, $\overline{D}_{\pi_m}+Z(\fra_{\xi,m})_{\pi_m}$ is effective.
	Thus $\overline{D}+{Z(\fra_{\xi,m})}\ge 0$ as a \textbf{b}-divisor, and 
	since ${Z(\fra_{\xi,m})}\le m D_\xi$, this implies that $\overline{D}+m D_{\xi}\ge 0$.
\end{proof}

\begin{definition}
	Let $v$ be a valuation on $X$, and $\xi=(X,v)$ the corresponding norm.
	We say that the valuation $v$ is \emph{\textbf{b}-divisorial} if $v(D)=\ord_{D_\xi}(D)$ for every divisor $D$ on $X$. Equivalently, $v$ is \textbf{b}-divisorial if $v(f)=\ord_{D_\xi}(\div(f))$ for every  $f\in K(X)$.
\end{definition}

\begin{proposition} \label{div-bdiv}
	If $v$ is a divisorial valuation, then it is \textbf{b}-divisorial.
\end{proposition}
\begin{proof}
	We can assume that $v=\ord_P$ where $P$ is a prime divisor in some model $X_\pi\overset{\pi}{\rightarrow} X$, and then $D_{\theta_P}=\Env_\pi(-P_\pi)$ as noted above.
	Therefore, the inequality of $\bbR$-\textbf{b}-divisors $\overline{D}+tD_{\theta_P}\ge 0$ implies $\pi^*(D)-tP_\pi \ge \pi^*(D)+t(D_{\theta_P})_\pi \ge 0$ and $v(D)\ge t$.
	Therefore $v(D)\ge \ord_{D_\xi}(D)$, and by \ref{pro:seminorm_bound_positivity} $v(D)= \ord_{D_\xi}(D)$ for every $D$. Thus $v$ is \textbf{b}-divisorial.
\end{proof}

\begin{proposition}

\label{pro:Dxi-cartier}
    $D_\xi$ is a nonzero Cartier \textbf{b}-divisor if and only if $\xi$ is a norm and $v_\xi$ is divisorial.
\end{proposition}

\begin{proof}
    If $\xi$ is a norm and $v_\xi$ is divisorial then $D_\xi$ is Cartier and nonzero as shown in Example \ref{exa:divisorial-envelope}.
    
    Let us now prove the converse, so assume $D_\xi$ is Cartier and nonzero.
    By Proposition \ref{pro:coefficients-dxi}, part \ref{part:strict_seminorm_zero}, since $D_\xi$ is nonzero, $\xi$ is a norm associated to a valuation $v=v_\xi$ of $K(X)$, and since $D_\xi$ is Cartier, it is meaningful to evaluate $v(D_\xi)$.
    We claim that $v(D_\xi)=-1$.
    If $D$ is an effective divisor on $X$, the inequality of $\bbR$-\textbf{b}-divisors $\overline{D}+tD_\xi\ge 0$ implies $v(D)\ge tv(-D_\xi)$. 
    Therefore $v(D)\ge \ord_{D_\xi}(D)v(-D_\xi)$, and since $\ord_{D_\xi}(D)\ge v(D)$ by Proposition \ref{pro:seminorm_bound_positivity}, we conclude that $v(D_\xi)\ge -1$, with an equality implying that $v = \ord_{D_\xi}$ (so $\ord_{D_\xi}$ is a valuation on $\cO_{X,\red(\xi)}$). 
    Let $X_\pi$ be a model where the Cartier divisor $D_\xi$ is determined, and write $T=(D_\xi)_{X_\pi}$.
    If $v(D_\xi)=-1$ then we have shown that $\ord_{D_\xi}$, the order of vanishing along $-T$, is a valuation on $\cO_{X,\red(\xi)}$; this implies that there is a prime component $P$ of $T$ with
    $\ord_{D_\xi}(\pi^*f)=\ord_P(\pi^*f)/\ord_{P}(-T)$ for every $f\in\cO_{X,\red(\xi)}$, and in particular $v$ is a divisorial valuation. 
        
     It remains to prove that indeed $v(D_\xi)=-1$, so fix a model $X_\pi$ where $D_\xi$ is determined, and let $T=(D_\xi)_{X\pi}$, $V=\red_{X_\pi}(\xi)$, and let $\xi_\pi$ the seminorm in $\frX_\pi$ determined by the valuation $v$ of $K(X_\pi)=K(X)$.
    By the birational nature of $D_\xi$ (\ref{pro:functoriality}) $D_\xi = \Env_\cX(D_{\xi_\pi})$, and in particular $D_{\xi_\pi}\ne 0$.
    
    If $V$ is a divisor we are done, so assume it is not and let $E_V$ be the exceptional divisor of blowing up $V$.
    We  know that $D_{\xi_\pi}$ is supported on divisors above $\red_{X_\pi}(\xi)$, so in particular $\alpha:=\ord_{E_V}(D_{\xi_\pi})< 0$, and by definition
    $$
    \alpha = \lim_{m\to\infty} \frac{\ord_{E_V} Z(\fra_{\xi_\pi,m})}{m} .
    $$
    The fact that $D_\xi$ is Cartier determined on $X_\pi$ means that 
    $$
    \lim_{m\to\infty}\frac{\ord_{E_V} Z(\fra_{\xi,m})}{m}=\ord_{E_V} D_\xi = \ord_{E_V} \overline{T}.
    $$
    Since $\fra_{\xi,m}=\pi_* \fra_{\xi_\pi,m}$, we get $Z(\fra_{\xi,m})\le Z(\fra_{\xi_\pi,m})$ and it follows that
    $$
    \lim_{m\to\infty}\frac{\ord_{E_V} (Z(\fra_{\xi_\pi,m})-m\overline{T})}{m} \ge
    \lim_{m\to\infty}\frac{\ord_{E_V} (Z(\fra_{\xi,m})-m\overline{T})}{m} =0
    $$
    Let $\pi':X_{\pi'}\rightarrow X_\pi$ be the normalized blowup of $\fra_{\xi_\pi,m}$. On that variety, we have
    $$
    \cO_{X_{\pi'}}(Z(\fra_{\xi_\pi,m})-\lfloor m \pi'^*(T)\rfloor)=
    (\fra_{\xi_\pi,m}: \cO_{X_\pi}(\lfloor mT\rfloor))\cdot \cO_{X_{\pi'}}
    $$
    and the definition of valuation ideals gives 
    $$
    \fra_{\xi_\pi,m}: \cO_{X_\pi}(\lfloor mT\rfloor) = \fra_{\xi_\pi, m+v (\lfloor mT \rfloor)}.
    $$
	Therefore
    \begin{gather*}       
    0\le\lim_{m\to\infty}\frac{\ord_{E_V} (Z(\fra_{\xi_\pi,m})-m\overline{T})}{m} =
    \lim_{m\to\infty}\frac{\ord_{E_V} Z(\fra_{\xi_\pi,m+v(\lfloor mD_\xi\rfloor)})}{m} = \\ \lim_{m\to\infty} \frac{Z(\fra_{\xi_\pi,m(1+v(D_\xi))})}{m}=
    (1+v(D_\xi))\lim_{m\to\infty} \frac{Z(\fra_{\xi_\pi,m})}{m}=
    (1+v(D_\xi)) \alpha.
    \end{gather*}
    Since $\alpha<0$ and $1+v(D_\xi)\ge 0$, it follows that $1+v(D_\xi)=0$ as claimed.
    \end{proof}

\begin{conjecture}\label{con:characterize-b-divisorial}
    $\xi$ is \textbf{b}-divisorial if and only if $D_\xi\ne 0$. 
\end{conjecture}
In all cases with $D_\xi\ne 0$ that we are able to compute, there is actually an equality $v_\xi=\ord_{D_\xi}$ (see section \ref{sec:continuity_Dxi} below).
Not long after this work was made available on arXiv, Lu Qi announced in \cite{qi26} a proof of Conjecture \ref{con:characterize-b-divisorial}.

\begin{example}
Let $X = \mathbb{A}^2$ with coordinates $(x,y)$ and maximal ideal $\mathfrak{m} = (x,y)$. It is not hard to construct a sequence of quasi-monomial valuations whose log discrepancies diverge.

\begin{itemize}
\item For $t \geq 1$, define $v_t$ by:
\begin{equation}
v_t(x) = 1, \quad v_t(y) = t,
\end{equation}
extended to $\mathcal{O}_{X,0}$ via:
\begin{equation}
v_t\left(\sum c_{a,b}x^a y^b\right) = \min\{a + t b \mid c_{a,b} \neq 0\}.
\end{equation}
These satisfy:
\begin{itemize}
\item $v_t(\mathfrak{m}) = \min(1,t) = 1$ for all $t$,
\item $A(v_t) = 1 + t \xrightarrow{t \to \infty} +\infty$,
\end{itemize}

\item The pointwise limit $v_\infty := \lim_{t \to \infty} v_t$ satisfies:
\begin{equation}
v_\infty(f) = \text{ord}_{x}(f) = \max\{k \mid x^k \text{ divides } f\}.
\end{equation}
This has:
\begin{itemize}
\item $A(v_\infty) = +\infty$ (non-quasi-monomial), but
\item $v_\infty(\mathfrak{m}) = 0 \neq 1$.
\end{itemize}
This shows that to have a well defined limit in the Berkovich space it is necessary to normalize the valuation. 
\item 
For example:
\begin{equation}
v_t' := \frac{v_t}{t}, \quad \text{i.e., } v_t'(x) = t^{-1}, \ v_t'(y) = 1.
\end{equation}
Now:
\begin{itemize}
\item $v_t'(\mathfrak{m}) = t^{-1} \xrightarrow{t \to \infty} 0$,
\item $A(v_t') = t^{-1}(1 + t) \xrightarrow{t \to \infty} 1$,

So that in this case $\lim_{t\to \infty} \frac{A(v_t')}{t}=0.$
\end{itemize}
\end{itemize}
\end{example}

\begin{proposition}\label{pro:sublinear-is-bdivisorial}
    Let $v$ be a valuation of $K(X)$ and $\xi=(X,v)$ the associated norm.
    Suppose that there exists a sequence of models $X_m$ and 
    regular systems of parameters $z^m_j\in \cO_{X_m,\red_{X_m}(v)}, j=1,\dots,c_m$
    \begin{enumerate}
        \item $Z(\fra_{\xi,m})$ and $\mathring Z_m(\fra_\bullet)$ are determined by
        SNC divisors on $X_m$ for every $m$,
        \item $-Z(\fra_{\xi,m})$ and $-\mathring Z_m(\fra_\bullet)$ are locally given by 
        monomial equations in $z^m$ at $\red_{X_m}(v)$,
        \item $\displaystyle\lim \frac{A(ev_{z^m}(v))}{m}=0$,
    \end{enumerate}
    where $ev_{z^m}:\frX\rightarrow (\bbR_{\ge 0}\cup\{+\infty\})^c$ is the retraction 
    associated to the system of parameters $z^m$ and the model $X_m$.
    
    Then $v=\ord_{D_\xi}$, i.e., $v$ is \textbf{b}-divisorial.
\end{proposition}

\begin{proof}
Let $p=\red_X(\xi)$ and let $(X_m, z^m)$ be a sequence of 
models and regular sequences satisfying the hypotheses.
Due to Proposition \ref{pro:seminorm_bound_positivity}, we only need to prove that
for every $f\in \cO_{X,p}$, $v(f)\ge \ord_{D_\xi}(f)$.
By Proposition \ref{pro:above_below_limit}, for every $m\ge 0$ the inequality
$\ord_{\mathring{Z}_m(\fra_{\xi,\bullet})}(f)\ge \ord_{D_\xi}(f)/m$ is satisfied,
and hence
$$
v(f)\ge \frac{\ord_{D_\xi}(f)}{m}\,v(-\mathring{Z}_m(\fra_{\xi,\bullet})).
$$
So the result will follow by proving that $\liminf v(-\mathring{Z}_m(\fra_{\xi,\bullet}))/m\ge1$.

Now we use the sequence of retractions $ev_{z^m}$; 
since for every $m$ the divisors $Z(\fra_{\xi,m})$ and $\mathring Z_m(\fra_\bullet)$ are 
determined by monomials on $X_m$, the definition of the retraction $ev_{z^m}(v)$ gives that
$$\frac{v(-\mathring{Z}_m(\fra_{\xi,\bullet}))}{m} = 
\frac{ev_{z^m}(v)(-\mathring{Z}_m(\fra_{\xi,\bullet}))}{m}$$ 
which due to \cite[6.2]{JM12} is bounded below by
$$\frac{ev_{z^m}(v)(-Z(\fra_{\xi,m}))-A(ev_{z^m}(v))}{m}=
\frac{v(-Z(\fra_{\xi,m}))-A(ev_{z^m}(v))}{m}=1-\frac{A(ev_{z^m}(v))}{m},
$$
and the claim follows.
\end{proof}

\begin{definition}\label{def:sublinear-discrepancy}
	The hypothesis in the preceding Proposition will be useful
	later on, hence we give a name to it.
	A valuation $v$ has \emph{sublinear log-discrepancy} if there
	is a sequence of models $X_m$ satisfying properties 1-3 in 
	Proposition \ref{pro:sublinear-is-bdivisorial}.
	Naturally, if $v$ has bounded log-discrepancy then it has sublinear log-discrepancy:
\end{definition}

\begin{theorem}
	If $A(v)<\infty$ then $v$ has sublinear log-discrepancy.
	In particular $v$ is \textbf{b}-divisorial.
\end{theorem}

\begin{proof}
	We only need to show that there exists a sequence of models $X_m$ and regular
	systems of paramters $z^m$ as in the previous proposition.
	Let $X_m$ be any model where $Z(\fra_{\xi,m})$ and $\mathring Z_m(\fra_\bullet)$ are determined by SNC Cartier divisors, and $z^m$ any regular system of parameters at $\red_{X_m}(v)$ such that $-Z(\fra_{\xi,m})$ and $-\mathring Z_m(\fra_\bullet)$ are
	given by monomials (which can be done as they are effective SNC divisors).
	By definition of the log discrepancy, $A(ev_{z^m}(v))\le A(v)$ for each model,
	so $\lim A(ev_{z^m}(v))/m=0$ as needed.
\end{proof}

\begin{corollary}
	If $v$ is an Abhyankar valuation, then it is \textbf{b}-divisorial.
\end{corollary}
\begin{proof}
    This can be proved directly using the fact proved by 
	Ein--Lazarsfeld--Smith in \cite{ELS03} that if $v=v_\xi$ is Abhyankar, then there exists $k$ (fixed for a given $v$) such that 
	$$\fra_{\xi, m}^\ell \subseteq \fra_{\xi, \ell m} \subseteq \fra_{\xi, m-k}^\ell$$
	for every $m\in \bbR_{\ge0}$ and every $\ell\in\bbN$.

    Alternatively, since by \cite[Proposition 5.1]{JM12}, Abhyankar valuations
    have finite log-discrepancy, it follows from the previous Theorem.
 \end{proof}

	Having sublinear log-discrepancy is a technical hypothesis which
simplifies the proofs of some statements which we conjecture
to hold for all \textbf{b}-divisorial valuations.

\begin{example}\label{exa:sublinear-multiplier}
	Let $v$ be a valuation with sublinear log-discrepancy,
	and let $\frj_m(\fra_{\bullet})$ be the subadditive sequence
	of asymptotic multiplier ideals.
	Then $\lim v(\frj_m(\fra_{\bullet}))/m =1$.
	If $v$ has bounded log-discrepancy, this was proved by Jonsson and Musta\c{t}\u{a} in \cite[6.2]{JM12}. 
	As seen in the proof of Proposition \ref{pro:sublinear-is-bdivisorial},
	their computation can be adapted to bound 
	$$\frac{v(\frj_m(\fra_{\bullet}))}{m} =\frac{v(-\mathring{Z}_m(\fra_{\xi,\bullet}))}{m} $$ 
	in the sublinear case as well, giving the desired limit.
\end{example}

\begin{lemma}\label{lem:supremum_toric}
	Let $\xi\in\frX$ be a seminorm, and let $\tilde P$ be a prime \textbf{b}-divisor. Then
	$$
	\ord_{\tilde P}(D_\xi)= \sup \{\ord_{\tilde P}(D_{ev_{z}(\xi)})\}
	$$
	where the supremum runs over all regular systems of parameters $z$ in all $\cO_{X_\pi,p}$ for all models $X_\pi$.
\end{lemma}
\begin{proof}
	First observe that if $\xi$ is not a norm, i.e., $\ker(\xi)$ is (the generic point of) a proper subvariety $V$, then for every smooth point $p_0$ of $V$ and every system of parameters $z$ with respect to which $V$ is locally monomial, $\ker(ev_{z}(\xi))=V$, so
	$$
	\ord_{\tilde P}(D_\xi)= 0 = D_{ev_{z}(\xi)} = \sup_{z} \{\ord_{\tilde P}(D_{ev_{z}(\xi)})\}
	$$
	by Proposition \ref{pro:coefficients-dxi}, part \ref{part:strict_seminorm_zero}.
	Thus in the rest of the proof we may assume that $\xi$ is a norm, corresponding to a valuation $v_\xi$ on $K(X)$.
	
	The inequality $\norm{f}_\xi\le\norm{f}_{ev_{z}(\xi)}$ holds for all $f\in \cO_{X,p}$ by \cite[Lemma 4.7]{JM12}, therefore $\fra_{\xi,m}\supseteq \fra_{{ev_{z}}(\xi),m}$ for every $m$ and so $D_{\xi} \ge D_{ev_{z}(\xi)}$. 
	Therefore $\ord_{\tilde P}(D_\xi)\ge \sup \{\ord_{\tilde P}(D_{ev_{z}(\xi)})\}$.
	
	We will show that a strict inequality 
	\begin{equation}
		\label{eq:nolimit}
		\ord_{\tilde P}(D_\xi)> \sup_{z} \{\ord_{\tilde P}(D_{ev_{z}(\xi)})\}
	\end{equation}
	would lead to a contradiction.
	By Proposition \ref{pro:coefficients-dxi},
	$$\ord_P(D_\xi)=-\inf_{f\in\cO_{X,\red(\xi)}}
	\frac{\ord_{P_\pi}(\pi^*(f))}{v_\xi(f)}=
	\sup_{f\in\cO_{X,\red(\xi)}}
	\frac{-\ord_{P_\pi}(\pi^*(f))}{v_\xi(f)}.
	$$
	Thus, the strict inequality \eqref{eq:nolimit} means that there exist $f\in \cO_{X,\red_{X}(\xi)}$ with 
	$$
	\frac{-\ord_{P_\pi}(\pi^*(f))}{v_\xi(f)} > \sup_{z} \{\ord_{\tilde P}(D_{ev_{z}(\xi)})\}
	$$
	Now, since $v_\xi(f)=\lim v_{ev_{z}(\xi)}(f)$, it follows that there exists some $z_0$ with 
	$$
	\frac{-\ord_{P_\pi}(\pi^*(f))}{v_{ev_{z_0}(\xi)}(f)} > \sup_{z} \{\ord_{\tilde P}(D_{ev_{z}(\xi)})\}=
	\sup_{z}\sup_{g\in\cO_{X,\red(\xi)}}
	\frac{-\ord_{P_\pi}(\pi^*(g))}{v_{ev_{z}(\xi)}(g)},
	$$
	a contradiction.
\end{proof}

\section{Continuity properties of the valuative divisor \texorpdfstring{$D_\xi$}{Dξ}}
\label{sec:continuity_Dxi}

\subsection{Birational nature of the valuative divisor \texorpdfstring{$D_\xi$}{Dξ}}
\label{sec:functoriality}

If $\phi:X\to Y$ is a generically finite morphism and $\xi\in\frX$, it is then natural to consider $\phi_*(\xi)$, the seminorm defined by pulling back functions, where $\norm{f}_{\phi_*\xi}=\norm{\phi^*(f)}_\xi$.

\begin{proposition} \label{pro:functoriality} Let $\phi:X\to Y$ be a birational morphism and let $\xi \in \frX$, then $$D_{\phi_*\xi} = \Env_{\cY}(\phi_*D_\xi).$$
\end{proposition}

\begin{remark}
Note that, since we are assuming that $\phi$ is birational, it induces a natural isomorphism of Riemann-Zariski spaces $\cX\cong\cY$, and 
\textbf{b}-divisors in $\Div \cY$ and $\Div \cX$ are also naturally identified via the isomorphism $\phi^*$ which ``forgets'' the traces of a \textbf{b}-divisor on all models that do not dominate $X$.
The inverse of $\phi^*$ consists in taking pushforward to those models, and so it coincides with the pushforward $\phi_*$ of \textbf{b}-divisors, defined in \cite[Section 1.5]{BdFF12} (see especially Proposition 1.15).  
Because of these identifications, in the sequel we will write $\Env_\cY(D_\xi):=\Env_{\cY}(\phi_*D_\xi)$ for brevity, as it can not possibly lead to any confusion.
The nef envelope operators $\Env_\cX$ and $\Env_\cY$ do differ however, because the conditions of being $X$-nef and $Y$-nef differ; thus, $D_\xi$ is $X$-nef by definition but need not be $Y$-nef, and so will generally differ from $\Env_\cY D_\xi$.
\end{remark}

\begin{proof}
    Let $\xi\in\frX$ be any seminorm. 
    By definition, $\ker \phi_*(\xi) = \phi (\ker \xi)$, so if $\xi$ is not a norm then $\phi_* \xi$ is not a norm either, both divisors $D_\xi$ and $D_{\phi_* \xi}$ are zero, and the equality holds trivially. 
    Thus we may assume that $\xi=(X,v_\xi)$ is the norm associated to a valuation $v_\xi$.

Note that $\Env_\cY(D_\xi)$ is well defined. Indeed, for every $m$, the \textbf{b}-divisor
$Z(\fra_{\xi,m})$ is Cartier, so $\Env_\cY(Z(\fra_{\xi,m}))$ is defined, is $X$-nef, 
and satisfies
$$
\Env_\cY(Z(\fra_{\xi,m}))\le Z(\fra_{\xi,m})\le D_\xi .
$$
Now for a given $m\in\bbR_{>0}$, we know by Corollary \ref{cor:Zxi-under-birational} that 
   $$\Env_{\cY}(Z(\fra_{\xi,m}))=\Env_{\phi\circ\pi_{m}}(Z(\fra_{\xi,m})_{X_m})=Z(\phi_*(\fra_m))=Z(\fra_{\phi_*(\xi),m})$$
   where $\pi_{m}:X_m\rightarrow X$ is the normalized blowup of $\fra_m$, and
   $$\mathring Z_m(\fra_{\phi_* \xi,\bullet}))=
   \mathring Z_m(\phi_*(\fra_{\xi,\bullet}))\ge\Env_\cY(\mathring Z_m(\fra_\bullet)).$$ 
Therefore by Proposition \ref{pro:above_below_limit} we have
$$
Z(\fra_{\phi_*(\xi),m})=\Env_{\cY}(Z(\fra_{\xi,m}))\le\Env_\cY(D_\xi)\le
\Env_\cY(\mathring Z_m(\fra_\bullet))
\le \mathring Z_m(\fra_{\phi_* \xi,\bullet})).
$$
Both the leftmost and rightmost terms in these inequalities converge to $D_{\phi_*\xi}$ as $m\to \infty$ so the claim follows.    
\end{proof}

\begin{remark}
    It should be possible, via an appropriate study of nef envelopes of push-forward Weil divisors, to extend Proposition \ref{pro:functoriality} to finite proper morphisms. 
    We do not need such a general \emph{functoriality} result in this work.
\end{remark}

\subsection{Global semicontinuity}

\begin{theorem} \label{thm:lowcont}
	The map $D:\frX \rightarrow \Div_\bbR(\cX)$ given by $\xi \mapsto D_\xi$ is lower semicontinuous.
\end{theorem}
\begin{proof}
    We only have to see that, for every prime \textbf{b}-divisor $\tilde P$, the map given by $\xi \mapsto \ord_{\tilde P}(D_\xi)$ is lower semicontinuous.
    So fix a prime \textbf{b}-divisor $\tilde P$ and a sequence or net $\{\xi_i\}_{i \in I}$ 
    coverging to $\xi=\lim_{i} \xi_i$; we need to show that 
    $\ord_{\tilde P}(D_\xi)\le \liminf_i \ord_{\tilde P}(D_{\xi_i})$.
    We argue by contradiction, assuming that     
    \begin{equation}\label{eq:no_semi_continuous_contradiction}
    0\ge e:=\ord_{\tilde P}(D_\xi) > \liminf_{i\in I} \ord_{\tilde P}(D_{\xi_i})=:\ell.    
    \end{equation}
    Let $\Delta = e-\ell>0$, and let $p=\red_X(\xi)$. 
    By Proposition \ref{pro:coefficients-dxi}, for every $e'<e$ 
    there exists $f\in \cO_{X,p}$ with 
    $$-
	    \frac{\ord_{\tilde P}(f)}{-\log\norm{f}_\xi}>e'.$$
    Choosing $e'=e-\Delta/4$ we get 
    $0\le\ord_{\tilde{P}}(f)<(e-\Delta/4)\log\norm{f}_\xi$ (note that in this case $\log\norm{f}_\xi<0$).
    Now since $\lim_i \xi_i = \xi$, we have that $\lim_i \log\norm{f}_{\xi_i}=\log\norm{f}_{\xi}$,
    so for every $\varepsilon>0$ 
    there exists $i_0\in I$ such that for $i\ge i_0$,
    the inequality $\log\norm{f}_{\xi_i}<\log\norm{f}_{\xi}(1-\varepsilon)<0$ holds. 
    But then, choosing $\varepsilon=\Delta/(2\Delta-4e)$ (which is positive since $e<0$ and $\Delta>0$), we have
    $$
    \ord_{\tilde{P}}(f)<\left(e-\frac{\Delta}{4}\right)\log\norm{f}_\xi<
    \left(e-\frac{\Delta}{4}\right)\frac{\log\norm{f}_{\xi_i}}{1-\epsilon}=
    \left(e-\frac{\Delta}{2}\right)\log\norm{f}_{\xi_i},
    $$
    which again by \ref{pro:coefficients-dxi} implies that $\ord_{\tilde P} D_{\xi_i}\ge e-\Delta/2=\ell+\Delta/2$ for every $i\ge i_0$,
    in contradiction with \eqref{eq:no_semi_continuous_contradiction}.
    This finishes the proof.
\end{proof}
\begin{remark}
    Note that Theorem \ref{thm:lowcont} applies in particular to cases with $D_\xi=0$.
    For them we have that, if $\lim_i \xi_i=\xi$ with $D_\xi=0$, then
    $$
    0=\ord_{\tilde P}(D_\xi) \le \liminf_{i \in I} \ord_{\tilde P}(D_{\xi_i}) \le 0
    $$
    for every $\tilde P$ by the anti-effectivity of all $D_{\xi_i}$.
    Therefore
    $
    \liminf_{i} \ord_{\tilde P}(D_{\xi_i}) = 0
    $
    which, again since $\ord_{\tilde P}(D_{\xi_i})\le 0$ for all $n$, means
    that the sequence (or net) $\ord_{\tilde P}(D_{\xi_i})$ actually converges to $0$.
    So the map $D$ is in fact \emph{continuous} at every $\xi$ with $D_\xi=0$.
\end{remark}

\subsection{Continuity on quasimonomial families}

Our next goal is to show that the map $D$ is in fact \emph{continuous} over large sections
of $\frX$, notably over the sets $\Delta_{\pi,p,z}$ of quasimonomial (semi)norms described 
in section \ref{sec:quasi-monomial}.
To this end, we will consider first the case of monomial seminorms, and then we will
apply the results of section \ref{sec:functoriality} about the dependency of the valuative 
divisor on the base model relatively to which it is defined to be nef. 

If $X_\pi\rightarrow X$ is some birational model above $X$, we will use the notation  	$D^{\pi}_\xi$ for the \textbf{b}-divisor associated to $\xi$ with respect to the model $X_{\pi}$, i.e.,  $D^{\pi}_\xi = \lim_{m\to\infty} Z(\fra_{\xi,m}^{X_\pi})$, where, as above, $\fra_{\xi,m}^{X_\pi}\subset \cO_{X_\pi}$ is the valuative ideal sheaf of $\xi$ on $X_\pi$.
Recall from section \ref{sec:quasi-monomial} that, if 
$p\in X_\pi$ is a point such that $\cO_{X_\pi,p}$ is regular of dimension $c$, and $z_1, \dots, z_c$ is a regular system of parameters, there is a subset $\Delta_{\pi,p,z}$ parametrizing all seminorms that are monomial with respect to $z$, which is homeomorphic to $(\bbR_{\ge 0}\cup \infty)^c$ and is a retraction of $\frX$.

\begin{proposition}
\label{pro:toric-continuity-model}
	Fix a model $X_\pi$ of $X$, a point $p\in X_\pi$ such that the local ring $\cO_{X_\pi,p}$ is regular of dimension $c\le \dim X$, and $z_j\in \cO_{X_\pi,p}$, $j=1,\dots, c$ a regular system of parameters.
	For every prime \textbf{b}-divisor $\widetilde{P}$, the function $(\bbR_{>0}\cup \{\infty\})^c \rightarrow \bbR$ defined by the composition of $w \mapsto \norm{\cdot}_w$ and $\xi \mapsto \ord_{\widetilde{P}}(D^{\pi}_\xi)$ is continuous in 
    $(\bbR_{>0}\cup \infty)^c\subsetneq \Delta_{\pi,p,z}$.
\end{proposition}
\begin{proof} This result is a consequence of \cite[Proposition 3.1]{JM12}. 
For every $m$ the ideal $\fra_{\xi_w,m}$ is monomial, so we only need to consider toric divisors $P$ at $p$ or on toric blowups above $p$ and the prime \textbf{b}-divisors $\widetilde{P}$ determined by them.

In particular, if all $w_i\ne 0$ this allows to bound the $P$-valuation of $\fra_{\xi_w,m}$ by the minimum of the toric pairings of $P$ with each ``monomial'' $z_i^{m/w_i}$, the actual value being computed by a monomial with integer exponents close to $z_i^{m/w_i}$.
Then the limit as $m\to\infty$ of $v_P(\fra_{\xi_w,m})/m$, which is the
opposite of the coefficient of $P$ in $D^{\pi}_{\xi_w}$, is the minimum of the formal toric pairings of $P$ with $z_i^{1/w_i}$.
This is clearly continuous as a function of $w\in (0,\infty]^c$ independently of $P$.
\end{proof}

Observe that the function above is \emph{never} continuous at any point of the boundary.
Indeed, let $I\subseteq \{1,\dots,c\}$ be the nonempty set of indices with $w_i=0$. 
If $I=\{1,\dots,c\}$ then $\xi_w$ is the trivial norm, $D^{\pi}_\xi=-\infty$, and clearly the map is not continuous at $w$, so we can assume that $\{1,\dots,c\}\setminus I$ is nonempty.
Let $P$ be the toric divisor associated to the blowup centered at the toric ideal $(z_i)_{i\not \in I}$; in this case the toric pairing is given by the ``indicatrix'' vector $\sum_{i\not\in I}\mathbf{e}_i$.
For every $i\in I$ and every $w'$ near $w$ with $w'_{i}\ne 0$, $z_i^{\lceil m/w_i \rceil} \in \fra_{\xi_w,m}$ has $\ord_{P}(z_i^{\lceil m/w_i \rceil})=0$, so $\ord_{\tilde P} D^{\pi}_{\xi_{w'}}=0$.
On the other hand, 
$$
\ord_{\tilde P}(D^{\pi}_{\xi_w}) = - \min_{i\not\in I} 1/w_i <0
$$
(since $w_i\ne 0$ for $i\not\in I$).
In other words, since at the limit, the ideals $\fra_{\xi_{w},m}$ contain no monomial involving only the $z_i, i\in I$, the toric pairing with $\mathbf{v}$ ``jumps'' to a strictly negative coefficient of $P$ in $D^{\pi}_{\xi_{w}}$.

Notice that if $c<n=\dim X$, then at each point $p'\in \overline{p}$ it is possible to extend the regular system of parameters $z=(z_1,\dots, z_c)$ to a regular system $z'=(z_1,\dots, z_c,\dots,n)$ for $\cO_{X_\pi,p'}$. 
Each such extension provides a larger set $\Delta_{\pi,p',z'}$ of seminorms, containing $\Delta_{\pi,p,z}$ on its boundary, and the collection of all such extensions cover all ``branches'' of the Berkovich space $\frX$ near $\Delta_{\pi,p,z}$.

\begin{example} \label{ex:jump} 
We will now describe an example of discontinuity and its connection with the non-continuity of the Nakayama--Zariski decomposition (also known as $\sigma$-decomposition).

Let us fix a model $X_\pi$ of $X$, a point $p \in X_\pi$ such that the local ring $\cO_{X_\pi,p}$ is regular of dimension $c \le \dim X$, and $z_j \in \cO_{X_\pi,p}$, $j = 1, \dots, c$ a regular system of parameters. Let $D_1 = \div_+(z_1), \dots, D_c = \div_+(z_c)$ be the divisors on $X_\pi$ associated to these parameters, and let $\Gamma_{D_1} = \ord_{D_1}$ be the divisorial valuation associated to the first divisor $D_1$, with $\theta_{D_1} = (X, \ord_{D_1}) \in \frX$ the associated norm, i.e., for $f \in K(X)$, 
\[
\norm{f}_{\theta_{D_1}} = \exp(-\ord_{D_1}(f)).
\]

Recall that, to consider the $\sigma$-decomposition, we need to work with effective divisors. Let us denote $\widetilde{D}_1$ the \textbf{b}-divisor obtained from $D_1$ via pushforward and strict transform, and $\overline{D}_1$ the Cartier closure of $D_1$. In particular, $\widetilde{D}_1 \leq \overline{D}_1$.

We now construct explicitly a sequence $\xi_k \in \frX$ such that:
\begin{itemize}
\item $\red(\xi_k) = p$ for all $k$,
\item In local coordinates, $\xi_k \in \bbR^c_{>0}$ and $\lim_{k\to \infty} \xi_k = \theta_{D_1}$,
\item $\overline{D}_1 + D_{\xi_k} \geq 0$ for all $k$.
\end{itemize}

For example, consider the sequence of quasimonomial valuations $\xi_k$ defined in coordinates $(z_1, \dots, z_c)$ by the weight vector $(1 + \frac{1}{k}, \frac{1}{k}, \dots, \frac{1}{k})$. 

Then:
\begin{itemize}
\item $\red(\xi_k) = p$ since all $\xi_k$ are centered at $p$,
\item $\xi_k \to \theta_{D_1}$ as $k \to \infty$ because the weights converge to $(1, 0, \dots, 0)$,
\item $\overline{D}_1 + D_{\xi_k} \geq 0$ since $D_{\xi_k}$ is effective and $\overline{D}_1$ is Cartier.
\end{itemize}

We obtain that $\Gamma_{D_1}(\overline{D}_1 + D_{\xi_k}) = 1$ for all $k$, but 
\[
\Gamma_{D_1}(\overline{D}_1 + D_{\Gamma_{D_1}}) = \Gamma_{D_1}(\overline{D}_1 - \widetilde{D}_1) = 0.
\]
\end{example}

\begin{corollary} \label{cor:delta-notation}
	The map $\bbR^c_{>0}\rightarrow \Div_\bbR(\cX)$ defined by the composition of $w \mapsto \norm{\cdot}_w$ and $\xi \mapsto D^{\pi}_\xi$ is continuous for the topology of coefficient-wise convergence on $\Div_\bbR(\cX)$.
\end{corollary}
\begin{proof}
	Follows from Proposition \ref{pro:toric-continuity-model}  by the definition of the topology on $\Div_\bbR(\cX)$.
\end{proof}

    Recall from section \ref{sec:quasi-monomial} that for every smooth model $X_\pi$, every point $p\in X_\pi$ and every regular system of parameters $z=(z_j)$ in $\mathcal{O}_{X_\pi,p}$ the map $$ev_{z}:\frX \rightarrow\Delta_{\pi,p,z}\cong(\bbR_{>0}\cup\{+\infty\})^c$$ defined in \eqref{eq:retraction} is a retraction  and $\frX$ is an inverse limit of such $\Delta_{\pi,p,z}$'s.

\begin{lemma} \label{pro:net}
       Let $W_i$ be a sequence (or net) of $\bbR$-Weil \textbf{b}-divisors which converges to an $\bbR$-Weil \textbf{b}-divisor $W$ coefficient-wise. Assume that each $W_i$ admits a nef envelope.
       Then $W$ admits a nef envelope which is bounded below by 
$$
W' := \liminf \Env_{\cX}(W_{i})
$$
\end{lemma}

\begin{proof}
Since for every $i$ we have an inequality $\Env_{\cX}(W_{i}) \le W_i$, the limit 
inferior satisfies $W'\le W$.
On the other hand, by closedness of the $X$-nef cone, $W'$ is $X$-nef. 
Thus $W$ admits a nef envelope, bounded below by $W'$.
\end{proof}

\begin{theorem}
\label{pro:toric-continuity}
	Fix a model $X_\pi$ of $X$, a point $p\in X_\pi$ such that the local ring $\cO_{X_\pi,p}$ is regular of dimension $c\le \dim X$, and $z_j\in \cO_{X_\pi,p}$, $j=1,\dots, c$ a regular system of parameters.
	For every prime \textbf{b}-divisor $\widetilde{P}$, the function $(\bbR_{>0}\cup \{\infty\})^c \rightarrow \bbR$ defined by the composition of $w \mapsto \norm{\cdot}_w$ and $\xi \mapsto \ord_{\widetilde{P}}(D_\xi)$ is continuous.
\end{theorem}
\begin{proof}
	By Theorem \ref{thm:lowcont}, the map is lower semicontinuous. We will show that 
	it is also upper semicontinuous.
	
    We saw in Proposition \ref{pro:toric-continuity-model} that the function
    $w \mapsto \ord_{\widetilde{P}}(D_{\norm{\cdot}_w^\pi})$ is continous, where $(D^{\pi}_\xi)$ denotes the \textbf{b}-divisor associated to $\xi$ with respect to the model $X_{\pi}$.

    By the birational nature of $D_\xi$ \ref{pro:functoriality}, $D_{\norm{\cdot}_w}=\Env_{\mathcal{X}}(D^{\pi}_{\norm{\cdot}_w})$ for every $w$, and 
    so the needed upper semicontinuity follows from Lemma \ref{pro:net}.
\end{proof}

\section{Positivity functions and their continuity}
\label{sec:positivity}
To define Seshadri constants for \textbf{b}-divisors, we rely on the intersection products from \cite{BdFF12} explained above.

Note that, by lemma \ref{lem:convergence}, the numerical class $[D_\xi]$ is the limit of all $[Z(\fra_{\xi,m})/m]$.
Indeed, by submultiplicativity all components of $Z(\fra_{\xi,m})_{X_{\pi}}$, $m>1$ are among the (finitely many) components of $Z(\fra_{\xi,1})_{X_{\pi}}$, for each model $X_{\pi}$ over $X$, so we can apply Lemma \ref{lem:convergence} to $W_j=Z(\fra_{\xi,j})/j$, with $V_i$ the span of the components of $Z(\fra_{\xi,1})_{X_{\pi}}$, to conclude that $[D_\xi]$ is the limit of $[Z(\fra_{\xi,m})/m]$.

Theorem \ref{thm:lowcont} implies the following:

\begin{corollary}
	The map $D:\frX \rightarrow N_{n-1}(\cX)$ given by $\xi \mapsto [D_\xi]$ is lower semicontinuous.
\end{corollary}
Note that the natural map $\Div_\bbR(\cX)\rightarrow N_{n-1}(\cX)$ (``taking numerical equivalence class'') is not continuous (see \cite[section 1.4]{BdFF12}) so this is not immediate from Theorem \ref{thm:lowcont}, but it needs Lemma \ref{lem:convergence}.
\begin{proof}
	Let $\xi_j$ be a sequence (or net) of seminorms in $\frX$ converging to a seminorm $\xi \in\frX$.
	We distinguish two cases. 
\begin{enumerate}
	\item The sequence (or net) $\xi_j$ has a subsequence (or net) consisting of seminorms \emph{divisorial on $X$}, i.e., such that the center $\red(\xi_j)$ is a divisor $D_j$ on $X$. 
	
	If the divisor $D_j$ is fixed, say $D_j=D$, for $j$ large enough, then $\norm{f}_{\xi_j}=\exp(-w_j)\ord_{D}(f)$ with $w_j\in[0,\infty]$ converging to $w\in[0,\infty]$, and $\norm{f}_{\xi}=\exp(-w)\ord_{D}(f)$. 
	In that case $[D_{\xi_j}]=-(1/w_j)[D]$ converge to $[D_\xi]=(-1/w)[D]$.
	
	Alternatively, there is no fixed divisor $D$ in $\cent(\xi_j)$ for $j\gg 0$.
	Therefore, $\lim \norm{D}_{\xi_j} = 0$ for every $D$ (because $\red(\xi_j)$ is a divisor $D_j$ and for almost all $j$, $D$ does not contain $D_j$). 
	Thus $\xi$ is the trivial norm, $D_\xi=-\infty$, and semicontinuity holds.
	
	\item For all $j$ large enough, the center $\red(\xi_j)$ is not divisorial on $X$.
	Then by Proposition \ref{pro:coefficients-dxi}, $(D_{\xi_j})_X=0$ for all $j$ large enough, and the result follows from Theorem \ref{thm:lowcont} and Lemma \ref{lem:convergence}.
	\qedhere
\end{enumerate}\end{proof}

\begin{definition}{\cite[Lemma 2.10]{BdFF12}}
$$ \Nef(\cX) = \lim_{\longleftarrow}{}_{\pi} \overline{\Mov(X_{\pi})}$$
where the limit is taken over all smooth (or $\mathbb{Q}$-factorial) models $X_{\pi}$, i.e. an
$\mathbb{R}$-Weil \textbf{b}-divisor $W$ is nef iff $W_{\pi}$ is movable on each smooth (or $\mathbb{Q}$-factorial) model $X_{\pi}$. 
\end{definition}

We leave to the reader to check the details of the goodness of the definition as explained in \cite{BdFF12}.
Next we introduce Seshadri constants in the context of \textbf{b}-divisorial valuations.

\begin{definition}[Seshadri constant and asymptotic order of vanishing of a \textbf{b}-divisorial valuation] \label{def:seshadri-D}
	Let $v=v_\xi$ be a \textbf{b}-divisorial valuation on $X$ and $D$ a nef Cartier \textbf{b}-divisor on $\cX$, then its Seshadri constant is
	
   $$\varepsilon(D, v):=\sup \left\{t \, |\, (D +t D_{\xi}) \text{ is nef }\right\}.$$
   
	Dually, we also define an asymptotic order of vanishing if $D$ is a big divisor:
   $$\omega(D, v):=\sup \left\{t \, |\, (D +t D_{\xi}) \text{ is pseudoeffective }\right\}.$$
\end{definition}
Here by pseudoeffectibe Weil \textbf{b}-divisor we mean a limit of effective Cartier \textbf{b}-divisors. 
By definition pseudoeffective \textbf{b}-divisors form a closed convex cone. 

The asymptotic order of vanishing \cite{DHKRS17, BKMS} is also called \emph{pseudoeffective threshold} in the literature, see e.g. \cite{BJ25}. Its reciprocal is commonly called \emph{Waldschmidt constant}, especially on projective spaces.

	For seminorms that are not \textbf{b}-divisorial we expect $D_\xi$ to vanish
(Conjecture \ref{con:characterize-b-divisorial}), so we suggest alternative definitions, explicitly constructing the constant as a limit.

\begin{definition}[Seshadri constant and asymptotic order of vanishing of a seminorm] \label{def:seshadri-limit}
Let $\xi \in \frX$ and $D$ a nef (resp. big) Cartier \textbf{b}-divisor on $\cX$, then

$$\varepsilon(D, \xi):=\limsup_{m}\left(\sup \left\{t \, |\, mD +t Z(\fra_{\xi, m}) \; {\rm is \, nef}\right\}  \right).$$

$$\omega(D, \xi):=\limsup_{m}\left(\sup \left\{t \, |\, mD +t Z(\fra_{\xi, m}) \; {\rm is \, pseudoeffective}\right\}  \right).$$
\end{definition}

Some remarks on these definitions are in order.

\begin{remark}
	Definitions \ref{def:seshadri-D} and \ref{def:seshadri-limit} mimic and extend the classical definition of Seshadri (resp. Waldschmidt) constant at a subvariety in terms of nefness (resp. bigness) on the blowup, which therefore coincide with our definition when the valuation or seminorm is divisorial.	
\end{remark}

\begin{remark}\label{rem:waldschmidt-big-eff} The asymptotic order of vanishing of a seminorm satisfies
\begin{align*}
	\omega(D, \xi)&=\limsup_{m}\left(\sup \left\{t \, |\, mD +t Z(\fra_{\xi, m}) \; {\rm is \, big}\right\}  \right)\\
	&=\limsup_{m}\left(\sup \left\{t \, |\, mD +t Z(\fra_{\xi, m}) \; {\rm is \, effective}\right\}  \right), 
\end{align*}
as in fact for every $m$ we have
\begin{align*}
\sup \left\{t \, |\, mD +t Z(\fra_{\xi, m}) \; {\rm is \, pseudoeffective}\right\}&=
\sup \left\{t \, |\, mD +t Z(\fra_{\xi, m}) \; {\rm is \, big}\right\}\\&=
\sup \left\{t \, |\, mD +t Z(\fra_{\xi, m}) \; {\rm is \, effective}\right\} .  
\end{align*}
\end{remark}

\begin{remark}
	We expect that for \textbf{b}-divisorial valuations, both definitions \ref{def:seshadri-D} and \ref{def:seshadri-limit} of the Seshadri constant should agree.
In Lemma \ref{lem:sublinear-seshadri-well-defined} below we prove the consistency of
both definitions  for valuations with sublinear log-discrepancy.
\end{remark}

\begin{lemma}\label{lem:nefness-for-xnef}
	Let $D'\ge D$ be two $\cX$-nef $\bbR$-Weil \textbf{b}-divisors with $(D')_X=D_X$.
	If $D$ is nef then $D'$ is nef.
\end{lemma}

\begin{proof}
	By closedness of the nef cones, it is not restrictive to assume that $D$ and $D'$ are $\bbQ$-Weil \textbf{b}-divisors, and taking a suitable multiple, that they are integral Weil
	\textbf{b}-divisors.
	
	Assume $D$ is nef, which by \cite[Lemma 2.10]{BdFF12} means that on every model
	$X_\pi$, $D_\pi$ is movable. We need to prove that the same holds for $D'$.
	Write $\Delta_\pi = D'_\pi-D_\pi$ which is an effective exceptional divisor.
	As $D_\pi$ is movable, every fixed component of $|D'_\pi|=|D_\pi+\Delta_\pi|$ 
	must be a component of $\Delta_\pi$.
    But $D'$ is $\cX$-nef, so $D'_\pi$ is $X$-movable, in particular it has no fixed 
	exceptional components.
\end{proof}

\begin{remark}
    In Definition \ref{def:seshadri-limit}, both limits superior are in fact limits. 
    For the asymptotic order of vanishing this is essentially \cite[2.6]{BKMS};
    let us give some details on the proof in the present setting.
    Set $t_m=\sup \left\{t \, |\, mD +t Z(\fra_{\xi, m}) \; {\rm is \, pseudoeffective}\right\}$ and consider the auxiliary number $a_m=m/t_m$ which is the least real number such that $a_mD+Z(\fra_{\xi, m})$ is pseudoeffective.
Then, for every positive $m$ and $n$ we have
$$
(a_m+a_n)D+Z_{m+n}\ge (a_m+a_n)D+Z_{m}+Z_n = (a_mD+Z_{m})+ (a_nD+Z_n)
$$
which is pseudoeffective, as the sum of two pseudoeffective divisors. 
Therefore $a_{m+n}\le a_m+a_n$, i.e., $a_m$ is subadditive and hence 
$$
\lim_{m\to\infty} t_m^{-1} = \lim_{m\to\infty}\frac{a_m}{m}=\inf\frac{a_m}{m}=\inf t_m^{-1}
$$
by Fekete's lemma.
Taking inverses, the claim follows. 

For the Seshadri constant, the same proof works, observing that the difference between the two sides in the inequality
$$
(a_m+a_n)D+Z_{m+n}\ge (a_m+a_n)D+Z_{m}+Z_n
$$
is exceptional, so Lemma \ref{lem:nefness-for-xnef} applies and nefness of the right hand side implies nefness for the left hand side.\qed
    \end{remark}

\begin{lemma}\label{lem:sublinear-seshadri-well-defined}
    For every \textbf{b}-divisorial valuation $v=v_\xi$, $\varepsilon(D, \xi)\le \varepsilon(D, v)$ and $\omega(D, \xi)\le\omega(D, v).$ 
    If $v$ has sublinear log-discrepancy, then $\varepsilon(D, \xi)=\varepsilon(D, v)$ and $\omega(D, \xi)=\omega(D, v).$
\end{lemma}
\begin{proof}

    The definition of $\varepsilon(D,\xi)$ implies that $D+\varepsilon(D,\xi)D_\xi$ is a limit of nef $\mathbb{R}$-Cartier \textbf{b}-divisors 
    $D+(t_i/{m_i})Z(\fra_{\xi, {m_i}})$ for a suitable sequence $m_i$
    with $\lim m_i=\infty$ and $\lim t_i=\varepsilon(D,\xi)$. 
    In particular, $D+\varepsilon(D,\xi)D_\xi$ is nef and 
    $\varepsilon(D, \xi)\le \varepsilon(D, v).$
    
    On the other hand, as explained in Example \ref{exa:sublinear-multiplier},
    $$\lim \frac{v(\frj_m(\fra_{\bullet}))}{m} =\lim \frac{v(-\mathring{Z}_m(\fra_{\xi,\bullet}))}{m} =1$$ 
    hence there exist rational numbers $c_m=\frac{a_m}{b_m}>1$ such that 
    $\lim_{m\to \infty} c_m=1$ and $\frj_m(\fra_{\bullet})^{a_m}\subset \fra_{\xi,m}^{b_m}$. 
    For these numbers it holds
    $$
    a_m \Env_{\cX}(\mathring{Z}_m(\fra_{\xi,\bullet})) < b_m Z(\fra_{\xi,m}), 
    $$
 	and since $\mathring{Z}_m(\fra_{\xi,\bullet})\ge m D_\xi$ and $D_\xi$ is $X$-nef,
 	we deduce that  
 	$$
    c_m D_\xi \le c_m \Env_{\cX}(\mathring{Z}_m(\fra_{\xi,\bullet}))/m < Z(\fra_{\xi,m})/m. 	
 	$$
 	
 	We claim that for every $t<\varepsilon(D,v)$ there exist $m\gg 0$ such that $mD+tZ(\fra_{\xi,m})$ is nef.
 	Since the $c_m$ converge to 1, for $m$ large enough the inequality $c_m\le \varepsilon(D,v)/t$ holds, and hence $D+tc_mD_\xi$ is nef.
 	Since $mD+tZ(\fra_{\xi,m})$ is $\cX$-nef and
	$$mD+tZ(\fra_{\xi,m})\ge mD+t mc_m(D_\xi)_\pi \ge
	m(D+tc_m D_\xi),$$
 	we can apply Lemma \ref{lem:nefness-for-xnef} 
 	and condlude that $mD+tZ(\fra_{\xi,m})$ is nef, as needed.

    The proof for the asymptotic order of vanishing is similar, but simpler, as Lemma \ref{lem:nefness-for-xnef} is not needed. We leave the details to the interested reader.
\end{proof}

\subsection{Semicontinuity of the asymptotic order of vanishing}

\begin{theorem} \label{thm:waldshmidt-continuity}
    For every big divisor $D$,
	the map $\omega: \frX \to \mathbb{R}$ defined by $\xi \mapsto \omega(D, \xi)$ is lower semicontinuous.
\end{theorem}

\begin{proof}
    Let $v_\xi = -\log\norm{\cdot}_\xi$. 
    If $\xi$ is a norm then $v_\xi$ is the associated valuation; if it is not a norm,
    then $v_\xi$ is merely a semivaluation taking values in $\bbR \cup \{\infty\}$.
		By Remark \ref{rem:waldschmidt-big-eff},
		\begin{gather*}
			\omega(D,\xi)=\sup\{m/k | H^0(X_\pi, \fra_{\xi,m}(\pi^*(kD)))\ne 0\} \\
			= \sup\{v_\xi(s)/k | s \in H^0(X_\pi, \cO_{X_\pi}(\pi^*(kD)))\}.
		\end{gather*}	
For every $s \in H^0(X_\pi, \cO_{X_\pi}(\pi^*(kD)))$, the map $\xi \mapsto v_\xi(s)/k$
is continuous.
Therefore $\omega(D,\cdot)$, being the supremum of continuous functions, is lower semicontinuous.
\end{proof}

\begin{proposition}\label{pro:semicont_on_cont}
    Let $\frY\subset\frX$ be any subset of
    \textbf{b}-divisorial valuations with sublinear log-discrepancy over which the map $\xi \mapsto D_\xi$ is continuous.
    Then the map $\omega:\frY \rightarrow \bbR$ defined by $\xi \mapsto \omega(D,\xi)$ is upper semicontinuous.

    In particular, for every model $X_\pi$ of $X$, every point $p\in X_\pi$ such 
    that the local ring $\cO_{X_\pi,p}$ is regular of dimension $c\le \dim X$, 
    and every regular system of parameters $z_j\in \cO_{X_\pi,p}$, $j=1,\dots, c$,
	the map $\omega: \operatorname{interior}(\Delta_{\pi,p,z}) \to \mathbb{R}$ defined by 
    $w \mapsto \omega(D, \xi)$ is continuous.
\end{proposition}
\begin{proof}
    Let us consider a sequence or net $\xi_i$ converging to $\xi$ in $\frY$.
    Since $D_\xi = \lim D_{\xi_i}$, we have
$$D+ \lim_{w_i\to w} \omega(D,\xi_i) D_{\xi} = \lim_{w_i\to w} (D + \omega(D,\xi_i) D_{\xi_i})$$
which is pseudoeffective (as the limit of pseudoeffective divisors), 
so $\omega(D, \xi)\geq \lim_i \omega(D,\xi_i)$.
\end{proof}

\subsection{Semicontinuity of the Seshadri constant}

    In the remaining part of the section we will partially prove semicontinuity of Seshadri constants with respect to the topology of the Berkovich space. Here we will give an example of why there cannot be a global semicontinuity of such constants and gives a strong motivation for the reason why in the literature the center of the considered valuations is generally fixed and most likely a point.

\begin{example}
\label{exa:non-semicontinuous-seshadri}
    Let us consider a monomial valuation on $\mathbb{P}^2$ centered at the intersection of two lines that we can locally think as $x=0, \; y=0$. We can then consider the sequence of monomial valuations of the form $$v_{w_m}(\sum a_{ij}x^iy^j)=\min\{i + j/m \, |\, a_{ij}\neq 0 \}.$$
    The limit of such valuations is clearly the divisorial valuation centered at the line $x=0$.

    On $\mathbb{P}^2$ the class of a line $D\sim l$ is ample and it is then natural to consider $\varepsilon(D, \xi_{w_m})$. 
    It is obvious to compute $\varepsilon(D, \xi_{l})=1$, while $\varepsilon(D, \xi_{w_m}) = 1/m \to 0$.
\end{example}

In this subsection all seminorms (or valuations) under consideration have bounded
log-discrepancy, and hence the distinction between the two definitions of Seshadri constants 
is irrelevant.

\begin{theorem}    
	\label{thm:toric-seshadri-continuity}
    Fix a model $X_\pi$ of $X$, a smooth closed point $p\in X_\pi$, and $z_j\in \cO_{X_\pi,p}$, $j=1,\dots, n$ a regular system of parameters.
    Then the map $\varepsilon: \operatorname{interior}(\Delta_{\pi,p,z}) \to \mathbb{R}$ defined by $\xi \mapsto \varepsilon(D, \xi)$ is continuous.
\end{theorem}

\begin{conjecture}
    We believe that Theorem \ref{thm:toric-seshadri-continuity} should also hold
    for $p$ any point such that the local ring $\cO_{X_\pi,p}$ is regular of dimension 
    $c\le \dim X$. However, our proof below using toric pairings only works if all relevant
    curves can be assumed to have image in $X_\pi$ meeting $p$ properly, which is not
    assured in the general case.
\end{conjecture}

\begin{proof}
Consider, for each given curve $C$ on $X_\pi$, the value
$$\varepsilon_C(D, \xi_w):= \limsup_m \frac{mD \cdot C}{- Z(\fra_{\xi_w, m})\cdot C} \in \bbR_\infty$$
where $Z(\fra_{\xi_w, m})\cdot C$ stands for the intersection product, on any model where $Z(\fra_{\xi_w, m})$ is defined, of this divisor with the strict transform of $C$ (by the projection formula, this intersection number is independent on the model chosen).
The obvious equality
$$\varepsilon_C(D, \xi_w)=\limsup_{m}\left(\sup_t \left(t \, |\, (mD +t Z(\fra_{\xi_w, m})) \cdot C \ge 0\right)  \right) $$
and the defintion of the Seshadri constant above show that 
\begin{equation*}
   \varepsilon(D, \xi_w)=
   \inf_C \varepsilon_C(D, \xi_w),
\end{equation*}
where the infimum is taken over all curves $C$ on $X_\pi$.
Note that curves in higher models whose image in $X_\pi$ is a point are not relevant for the computation of $\varepsilon$, because $Z(\fra_{\xi_w, m})$ is $X_\pi$-nef.

We next analyze the dependence of $\varepsilon_C(D, \xi_w)$ on $w$, 
and we shall prove that there exist
a neighborhood $U_w$ and a constant $k_w$ such that for every curve $C$ relevant on $U_w$ for the computation of the Seshadri constant, the map 
$w\mapsto \varepsilon_C(D, \xi_w)$
is $k_w$-Lipschitz continuous on $U_w$. 
This will imply that $\varepsilon(D,\xi)$ is also $k_w$-Lipschitz continuous on $U_w$, i.e., it is locally Lipshitz continuous.

To make this argument precise we introduce a technical definition:
\begin{definition}
    Let $U\subset \Delta_{\pi,p,z}$.
    We say that a collection of curves $\mathcal{C}$ on $X_\pi$ is Seshadri $U$-sufficient if $\varepsilon(D, \xi_w)=   \inf_{C\in \mathcal{C}} \varepsilon_C(D, \xi)$ for every $w\in U$.
\end{definition}

The statement we need to prove is then the following:
\begin{proposition}\label{pro:cont-epsilon}
   There exists a neighborhood $U_w$, a constant $k_w$, and a Seshadri $U_w$-sufficient collection of curves $\mathcal{C}_w$ such that for every curve $C\in \mathcal{C}_w$, the map
   $$
   w \mapsto \varepsilon_C(D, \xi_w)
   $$
 is $k_w$-Lipschitz continuous.
\end{proposition}

As soon as Proposition \ref{pro:cont-epsilon} is proved, since the infimum of $k$-Lipschitz continuous functions is $k$-Lipschitz continuous the claimed continuity will follow.
\end{proof}

The proof of Proposition \ref{pro:cont-epsilon} is quite simpler in the case $X=X_\pi$; we prove first this case as Lemma \ref{lem:cont-Zm-above}. The difference with the general case is due to the presence of components of the divisors $Z_{\fra_{\xi_w,m}}$ contracted by $\pi$, and the continuity of these components' coefficients is dealt with in Lemma \ref{lem:cont-Zm-all}. With the help of these preliminary lemmas the proof of \ref{pro:cont-epsilon} will be significantly streamlined.

\begin{lemma}\label{lem:cont-Zm-above}
    For every $w\in \operatorname{interior}(\Delta_{\pi,p, z})$ there exists a neighborhood $U_w$, a constant $k_w$, and a Seshadri $U_w$-sufficient collection of curves $\mathcal{C}_w$ such that for every curve $C\in \mathcal{C}_w$, the map
   $$
   w' \mapsto \varepsilon_C^{X_\pi}(D, \xi_{w'})
   $$
   is is $k_w$-Lipschitz continuous on $U_w$.
\end{lemma}
\begin{lemma}\label{lem:cont-Zm-all}
    For every $w\in \operatorname{interior}(\Delta_{\pi,p, z})$ there exist a neighborhood $U_w$ and a constant $k_w$ such that for every prime divisor $E$ on $X_\pi$ contracted by $\pi$, the coefficient of $E$ in $D_{\xi_{w'}}$ is $k_w$-Lipschitz continuous as a function of $w'\in U_w$.
\end{lemma}

\begin{proof}[Proof of Lemma \ref{lem:cont-Zm-above}]
For the purposes of this lemma, all Seshadri constants and $Z$-divisors refer to the base model $X_\pi$, so we simply write $\varepsilon(D, \xi)=\varepsilon^{X_\pi}(D, \xi)$ and $Z(\fra)=Z^{X_\pi}(\fra)$.
    
Note first that if $C$ is a curve on $X_\pi$ which does not pass through $p$ $Z\left(\fra_{\xi_w,m}\right)\cdot C= 0$ for every $w$ and $m$, so $\varepsilon_C(D,\xi_w)$ is infinite and we may disregard it for the computation of Seshadri constants on a neighborhood of $w$ (i.e., $C$ will not be included in our Seshadri $U_w$-sufficient collection).
Now let $C$ be any curve that is potentially relevant, i.e., assume that $C$ has $k>0$ branches through $p$. 
The $k$ branches can be parameterised as $(f_1^{(i)}(t), \dots, f_n^{(i)}(t))$, where each $f$ is a holomorphic function on a disc $t\in U\subset \bbC$ around 0. Let $a_j = \sum_{i} \ord_t(f_j^{(i)})$ for $j=1,\dots,n$ and $\mathbf{a}=(a_1,\dots,a_c)$.

Next observe as in the proof of \ref{pro:toric-continuity-model} that, for every $m$, the ideal $\fra_{\xi_w,m}$ is monomial, generated by those monomials $z^\alpha$ with $\langle w,\alpha\rangle \ge m$.
Then, $Z(\fra_{\xi,m})$ is supported on toric divisors $P$ on toric blowups above $p$, and its intersection with $C$ is in fact equal to the opposite of the colength of $\fra_{\xi,m}|_C$, which can be computed torically as
$$
-\min\{\langle \mathbf{a},\alpha\rangle \,|\, \langle w,\alpha\rangle \ge m\}.
$$
For large $m$, this minimum will be attained by a monomial $x^\alpha$ with
$\alpha/m$ close to $1/w_j \mathbf{e}_j$ for some $j\in\{1,\dots,c\}$ (where we use $\mathbf{e}_j$ for the $j$-th basis vector in $\bbR^c$).
In any event,
$$
\varepsilon_C(D,\xi_w)=\frac{D\cdot C}{\min_{j=1,\dots,c} a_j/w_j } = 
D\cdot C \max_{j=1,\dots,c} w_j/a_j\, .
$$

Now fix a curve $C_0$ through $p$ and such that $\max\{w_j/a_j(C_0)\}$ is attained for a single $j=j_0$.
Then 
$$
\varepsilon_{C_0}(D,\xi_{w'}) = \frac{D\cdot C_0}{a_{j_0}(C_0)}\cdot {w'_{j_0}}
$$
for every $w'$ in a neighborhood $U_w$ of $w$.

Since $w$ belongs to the interior of $\Delta_{\pi,p,z}$, $w_i> 0\,\forall i$, so we may pick an arbitrary neighborhood $U_w$ and there exist constants $M_1,\dots,M_c>0$ such that $w'_i \ge M_i w'_{j_0}$ for every $w'\in U$ and every $i$.
Then, the following set of curves is Seshadri $U_w$-sufficient:
\begin{equation*}
\left\{ C \left| 
    \frac{D\cdot C}{a_j(C)}w'_j
   \le \frac{D\cdot C_0}{a_{j_0}(C_0)}w'_{j_0} \forall j
   \text{ s.t.} \frac{w'_j}{a_j(C)}= \max_i \frac{w'_i}{a_i(C)}
   \right.\right\}.
\end{equation*}
We can enlarge the set somewhat (so that it stays Seshadri $U_w$-sufficient) using a condition which is independent on $w'\in U_w$:
\begin{equation*}
  \mathcal{C} := \left\{ C \left| 
   \frac{D\cdot C}{a_j(C)}
   \le \frac{D\cdot C_0}{M_j a_{j_0}(C_0)} \forall j
   \right.\right\}.
\end{equation*}
In this way it is clear that that the function
$$
\varepsilon_{C}(D,\xi_w) = D\cdot C \max_j \frac{w'_j}{a_{j}(C)}
$$
is $k_w$-Lipschitz continuous on $U_w$ where  $k_w=\frac{D\cdot C_0}{\min(M_j)a_{j_0}(C_0)}$ for every $C\in \mathcal{C}$. 
\end{proof}

\begin{proof}[Proof of Lemma \ref{lem:cont-Zm-all}]
Given two vectors $w,w'$ in the interior of $\Delta$, define 
$$    
\alpha(w,w')=\max_{i=1,\dots,c} {w_i}/{w'_i}.
$$
The description above of the monomial ideals $\fra_{\xi_w,m}^{X_\pi}$ shows that, for every $m$, the following inclusions hold:
$$ \fra_{\xi_w,m\alpha(w,w')}^{X_\pi} \subseteq \fra_{\xi_{w'},m}^{X_\pi},   \qquad
\text{ and hence }
\fra_{\xi_w,m\alpha(w,w')} \subseteq \fra_{\xi_{w'},m}.
$$
These in turn imply inequalities
$$Z_{\xi_w,m\alpha(w,w')}\le Z_{\xi_{w'},m}$$
valid for all $m$, and it follows that 
\begin{equation}
    \label{eq:alpha-Dxi}
    \alpha(w,w')D_{\xi_{w}} \le D_{\xi_{w'}}\, .
\end{equation}

For a given $w$ in the interior of $\Delta$, choose a neighborhood $U_w$ and constants $A,B>0$ such that for every $w'=(w'_1,\dots,w'_n)\in U_w$, $A<w'_i<B$ for every $i$.

Again, the description above shows that $(z_1,\dots, z_n)^{\lceil m/A\rceil} \subset \fra_{\xi_w,m}$ for every $w \in U_w$ (where $\lceil x \rceil$ stands for the least integer greater or equal to $x$).

Fix $m_0\ge A/B$
and let $Z=Z((z_1,\dots, z_n)^{\lceil m_0/A\rceil})_{X_\pi}$,
so that for every $w'\in U_w$, $Z\le Z_{\xi_{w'},m_0}$ and hence
$Z\le m_0D_{\xi_{w'}}$.
For every $a, b$ positive real numbers, there are inequalities 
\begin{equation}
\label{eq:izumi}
 aZ + Z_{\xi_{w'},b} \le Z_{\xi_{w'}, am_0+b}   \, .
\end{equation}
and therefore
\begin{equation}
\label{eq:izumi-Dxi}
 aZ + bD_{\xi_{w'}} \le (am_0+b)D_{\xi_{w'}}   \, .
\end{equation}

Now let $\gamma(w',w'')=\max(0,\alpha(w',w'')-1)$, and note that
$\gamma(w',w'')\le \norm{w'-w''}_1/A$ in $U_w$. For every $w',w''\in U_w$ we have
$$
D_{\xi''}\ge \alpha(w',w'')D_{\xi_{w'}}\ge 
D_{\xi_{w'}}+\frac{\gamma(w',w'')}{m_0}Z \ge 
D_{\xi_{w'},m}+\frac{\norm{w'-w''}_1}{Am_0}Z 
$$

Symmetrically, we also have $D_{\xi_{w'}}\ge
D_{\xi_{w''},m}+\frac{\norm{w'-w''}_1}{Am_0}Z $ and so,
all components of $Z$ (which are the components of the $D_{\xi_{w'}}$ contracted by $\pi$) appear in $D_{\xi_{w'}}$ with coefficients that are Lipschitz continuous functions of $w'\in U_w$ with Lipschitz constant at most $1/Am_0$ times the largest coefficient of $Z$.
\end{proof}

\begin{proof}[Proof of Proposition \ref{pro:cont-epsilon}]
Fix, as in lemma \ref{lem:cont-Zm-all}, a contracted divisor \\
$Z=Z((z_1,\dots, z_n)^{\lceil m_0/A\rceil})_{X_\pi}$.
Let $C$ be any curve on $X_\pi$ not contained in $|Z|$, the support of $Z$.
$C$ has $k\ge 0$ branches through $p$.
If $k=0$, i.e., $C$ does not go trhough $p$, then $Z\left(\fra_{\xi_w,m}\right)\cdot C= 0$ for every $w$ and $m$, so $\varepsilon_C(D,\xi_w)$ is infinite and we may disregard it for the computation of Seshadri constants on a neighborhood of $w$ (i.e., $C$ will not be included in our Seshadri $U_w$-sufficient collection).

As in the proof of lemma \ref{lem:cont-Zm-above}, the intersection of $C$ with $Z_{\xi_w, m}$ can be computed on $X_\pi$ as the opposite of the colength of $(\fra_{\xi_w,m} \cdot \cO_{X_\pi})|_C$.
Now, the ideal $(\fra_{\xi_w,m} \cdot \cO_{X_\pi})$ is the intersection of a divisorial part supported on $|Z|$ and the monomial ideal described in \ref{lem:cont-Zm-above}; equivalently, it is the product of the divisorial part and the monomial part $\fra_{\xi_w,m_w'}^{X_\pi}$ where $m-m_w'$ is the $\xi_w$-valuation of the divisorial part. 
When $m$ tends to $\infty$, $\lim \frac{m-m'_w}{m} = v_{\xi_w}((D_{\xi_w})X_{\pi})$. 
Therefore,
$$
C \cdot D_{\xi_w} = -\min_{j=1,\dots,c} \frac{a_j}{w_j}
$$

The colength of the restricted ideal can be computed as the sum of the colengths of the factors, which are continuous as functions of $w$, and we are done.
\end{proof}

Lower semicontinuity extends to larger subsets of the space of seminorms, due to the following proposition, but it is not easy to characterize such subsets.

\begin{proposition}    
\label{pro:seshadri-supremum-models}
Let $\xi\in\frX$ be a seminorm, and let $D$ be a nef Cartier \textbf{b}-divisor. Then
$$
\varepsilon(D,\xi)= \sup \{\varepsilon(D,ev_{z}(\xi))\}
$$
where the supremum runs over all regular systems of parameters $z$ in all $\cO_{X_\pi,p}$ for all models $X_\pi$.

In particular, if $\xi\in \frX$ is such that $ev_{z}(\xi)$ belongs to the interior of
$\Delta_{\pi,p,z}$ for every $z$, then $\varepsilon(D,\cdot)$ is lower semicontinuous at $\xi$.
\end{proposition}

Note that there do exist norms $\xi$ satisfying the latter condition, and that the statement for such $\xi$ follows as there $\varepsilon(D,\cdot)$ is the supremum of continuous functions.

\begin{proof}
	By Lemma \ref{lem:supremum_toric}, $D_\xi$ is the limit of $D_{ev_z(\xi)}$.
	For every $t$ such that $mD+tZ(\fra_{ev_z(\xi),m})$ is nef, by Lemma \ref{lem:nefness-for-xnef}, $mD+tZ(\fra_{\xi,m})$ is nef as well.
	Therefore $	\varepsilon(D,\xi)\ge \sup \{\varepsilon(D,ev_{z}(\xi))\}$.
	However, a strict inequality leads to contradiction, by the same argument
	as in lemma \ref{lem:sublinear-seshadri-well-defined}.
\end{proof}

\section{The case of Surfaces} 
\label{sec:surfaces}
In this subsection, we determine $D_\xi$ for every $\xi$ when $X$ is a smooth surface.
Because of Proposition \ref{pro:coefficients-dxi},\ref{part:strict_seminorm_zero}, we need only consider the case of a real valuation $v$ on $K(X)$ and its associated norm $\xi=(X,v)$. 
If the center $P\subset X$ of $v$ on $X$ is divisorial, then $v=a\cdot \ord_P$ for some $a>0$ and one immediately checks that $D_\xi=-(1/a)\overline{P}$; so we assume in the rest of the section that the center of $v$ in $X$ is a point $p=\red_X(\xi)$.

We recall, following \cite{Cas00}, the classification of valuations centered at a point of a smooth surface using its sequence of centers (see also \cite{Spi90}).
The valuation $v$ determines a sequence of point blowing ups
\[
\dots \overset{b_{p_i}}{\longrightarrow} 
X_{p_{i-1}}\overset{b_{p_{i-1}}}{\longrightarrow} \dots
\overset{b_{p_2}}{\longrightarrow} X_{p_1}
\overset{b_{p_1}}{\longrightarrow} X
\]
where $p_1=p$ and in general $p_i=\red_{X_{p_{i-1}}}(\xi)$ is the center of $v$ on $X_{p_{i-1}}$; $b_q$ denotes the blowing up centered at $q$.
If $v$ is a divisorial valuation, then eventually $\red_{X_{p_{i}}}$ is a divisor, so the sequence of centers and blowing ups is finite; otherwise it is infinite.
As customary, we use the notation $E_i$ for the exceptional divisor of the $i$-th blowing up $b_{p_i}$ and $\tilde{E}_i$ for the corresponding prime \textbf{b}-divisor. 
We also fix the notation $\pi_i:X_{p_i}\rightarrow X$ for the composition of the first $i$ blowing ups.
There is a natural relation, called \emph{proximity}, between the centers: $p_i$ is proximate to $p_{j}$ if $p_i\in (\tilde{E}_{j})_{X_{p_{i-1}}}$.
Every center $p_i$, $i>1$, is clearly proximate to its predecessor $p_{i-1}$, and it can be proximate to at most another $p_{j}$ with $j<i-1$; in that case, $p_i$ is said to be a \emph{satellite}, otherwise it is \emph{free}.

There exist valuations for which there is a center $p_j$ with infinitely many centers $p_i, i>j$ proximate to it; however these are not real valuations by \cite[8.1.14]{Cas00}, so we assume henceforth that the set of $p_i$ proximate to each given $p_j$ is finite.
Consider the numbers $v_i=v(E_i)$. They are determined by $v_1$ and the so-called \emph{proximity equalities} \cite[8.1.7]{Cas00}
\begin{equation}\label{eq:proximity}
    v_j=\sum_{p_i\text{ prox. to }p_j} v_i.
\end{equation}
The fact that $v$ is a real valuation similarly implies that, for each given divisor $D$, only finitely many centers belong to the strict transform of $D$, and hence $v(D)$ can be evaluated using Noether's equality \ref{lem:noether}:
$$v(D)=\sum_{i\ge1} v_i \mult_{p_i}(\tilde{D})_{X_{p_{i-1}}},$$
in particular $v$ is determined by $v_1$ and the sequence of centers.

Now we can describe explicitly the ideals $\fra_{\xi,v}$ and hence the divisor $D_\xi$.

\begin{lemma}\label{lem:base-centers}
    Let $v$ be a real valuation on $X$ and $\xi=(X,v)$ the associated norm.
    For every $m\in \bbR$, there is an $n$ such that $X_{p_n}$ is the minimal resolution of $X_{\pi_m}$. 
\end{lemma}

\begin{proof}
    Since $\fra_{\xi,m}$ is a complete ideal, it is determined by a (finite) cluster of infinitely near points, based at $p_1$, and the blowup $X_{\pi_m}$ is a sandwiched singularity whose minimal resolution is obtained from $X$ by blowing up all points in the cluster of base points of $\fra_{\xi,m}$ \cite{Lip69,FS03}. 
    We proceed by induction on the number $n$ of base points of $\fra_{\xi,m}$.

    Let $e_1=\ord_{p_1}(\fra_{\xi,m})$ be the minimal multiplicity at $p_1$ of an element in $\fra_{\xi,m}$. 

    The ideal $\fra_{\xi,m}$ has a unique base point (i.e., $n=1$, in which case the claim is obvious) exactly when it is a power of the maximal ideal, $\fra_{\xi,m}=\mathfrak{m}_{p_1}^{e_1}$.
    This happens if and only if every element in $\cO_{X_{p_1}}$ with multiplicity at least $e_1$ belongs to $\fra_{\xi,m}$, i.e., $e_1 v_1 \ge m$.    

        So assume $e_1v_1< m$. There exist elements $f$ in $\fra_{\xi,m}$ with multiplicity exactly $e_1$, and in this case, by Noether's equality \ref{lem:noether}, $b_{p_1}^*(f) - e_1 E_1$ has positive valuation, i.e., it goes through $p_2$, for every such $f$.
    More precisely, if we let $\fra'_{\xi,x}$ be the ideal in $\cO_{X_{p_1},p_2}$ formed by functions with valuation at least $x$, Noether's equality gives
    \begin{equation}\label{eq:fra-one-blowup-step}
      \fra_{\xi,m}=(b_{p_1})_*(\fra'_{\xi,m-e_1v_1}(-e_1E_{p_1})).  
    \end{equation}
    It follows that the base points of $\fra_{\xi,m}$ are exactly $p_1$ and the base points of $\fra'_{\xi,m-e_1v_1}$.
    Now by the induction hypothesis, these base points are centers of $v_\xi$, and $X_{p_n}$ is the resolution of the normalized blowup of $X_{p_1}$ centered at $\fra'_{\xi,m-e_1v_1}$.
    Moreover, the equality \eqref{eq:fra-one-blowup-step} tells us that the normalized blowup of $X_{p_1}$ centered at $\fra'_{\xi,m-e_1v_1}$ factors through $X_{\pi_m}$
    (and they coincide if $E_{p_1}$ is not contracted in $X_{\pi_m}$). 
\end{proof} 

Next we determine the decomposition in prime divisors of $D_\xi$ for a divisorial valuation $v$. 
In this case there is a finite number of centers, which we denote by $n$, and for the last 0-dimensional center $p_n$, it holds $v=v_n\cdot \ord_{E_{p_n}}$.

\begin{lemma}\label{lem:am-divisorial}
    Let $v$ be a divisorial valuation with $n$ 0-dimensional centers 
    with values $v_1$, \dots, $v_n$, and let $\xi=(X,v)$ be the associated norm.
    Then for every $m\in \bbR$ such that $mv_n/\vsq$ is an integer,
    we have:
    \begin{gather*}
     Z(\fra_{\xi,m})_{\pi_n} = -\frac{m}{\vsq}\sum_{i\le n} v_i \overline{E_i}_{X_{p_n}}.  
    \end{gather*}
    where $\vsq=\sum v_i^2$.
\end{lemma}

Note that certainly there exist $m\in\bbR$ such that $mv_n/\vsq$ is an integer, and due to \ref{eq:proximity}, this implies that for every $i=1,\dots,n$, $mv_i/\vsq$ is an integer.
\begin{proof}
    Since $v=v_n\cdot \ord_{E_{p_n}}$, by definition we have
    $$
    \fra_{\xi,m}=(\pi_n)_*\left(\cO_{X_{p_n}}\left(-\frac{m}{v_n}E_n\right)\right).
    $$
    Let $D$ be the divisor $-\sum_{i\le n} \frac{mv_i}{\vsq} \overline{E_i}_{X_{p_n}}$ and $\fra_m=(\pi_n)_*(\cO_{X_{p_n}}(D))$.
    We need to show that $\fra_{\xi,m}=\fra_m$.
    
    It is immediate to check that $D\cdot (\tilde E_i)_{X_n}=0$ for every $i<n$. 
    By \cite[Lemma 3]{CS93} this implies that 
    $$
    -\ord_{E_n}(D)\cdot E_n = (-\ord_{E_n}(D)\cdot E_n-D) + D
    $$
    is the relative Zariski decomposition of $-\ord_{E_n}(D)\cdot E_n$ and therefore
    $$
    \fra_m = (\pi_n)_*(\cO_{X_{p_m}}(\ord_{E_n}(D)\cdot E_n)).
    $$
    On the other hand, since $\ord_{E_n}(\overline{E_i})=v_i/v_n$, one sees that 
    $$\ord_{E_{p_n}}(D)=-\sum \frac{mv_i}{\vsq}\frac{v_i}{v_n}=-\frac{m}{v_n}
    $$
    so $\fra_m=\fra_{\xi,m}$ and $D=Z(\fra_{\xi,m})$ as claimed.
\end{proof}
\begin{proposition}
\label{cor:divisorial-D}
     Let $v$ be a divisorial valuation  and let $\xi=(X,v)$ be the associated norm.
     Then $\vol_X(v)=1/\vsq$ and $$D_\xi=-{\vol_X (v)}\sum v_i \overline{E_i}.$$
\end{proposition}
Later we will show that this formula holds in all generality.

Recall from \cite{ELS03} that the volume of $v$ is 
$$\vol_X(v)=\lim_{m\to\infty} \frac{\length\cO_{X,p}/\fra_{m}}{m^2/2}.$$    

\begin{proof}
By the previous lemma, for every $m$ such that
$mv_n/\vsq$ is an integer the Hoskin--Deligne--Casas-Alvero formula \cite[6.1]{Cas00}gives
$$
\length\cO_{X,p}/\fra_{m}= \sum_{i=1}^n \frac{\frac{mv_i}{\vsq}(\frac{mv_i}{\vsq}+1)}{2}
=\frac{m^2}{2\vsq} + \frac{m\sum v_i}{2(\vsq)^2}.
$$
Dividing by $m^2/2$ and taking the limit we get $\vol_X(v)=1/\vsq$ as claimed.
On the other hand the previous lemma also gives
$$
D_\xi = \lim_{m\to\infty} Z(\fra_{\xi,m})/m=
-\frac{1}{\vsq}\sum_{i\le n} v_i \overline{E_i}_{X_{p_n}}=
-\vol(v)\sum_{i\le n} v_i \overline{E_i}_{X_{p_n}}.
$$
\end{proof}

\begin{theorem} \label{thm:b-divisorial surfaces}
     Let $v$ be a valuation on $X$ with zero-dimensional center, and let $\xi=(X,v)$ be the associated norm.
    Then $$D_\xi=-{vol_X (v)}\sum v_i \overline{E_i},$$
    where $v_1, v_2, \dots$ are the values of $v$ at its sequence of centers.
\end{theorem}

\begin{proof}
    Because of the Proposition above, we need only consider the case that $v$ is not divisorial. Our proof proceeds in two steps:
    \begin{enumerate}
        \item \label{step_sequence} We construct a sequence of divisorial valuations $v^{(k)}$ which is nondecreasing (i.e., $v^{(k)}(f)\le v^{(k+1)}(f)$ for every $f\in \cO_{X,p_1}$ and every $k\ge 1$), converges to $v$, and satisfies $$\lim_{k\to \infty}D_{\xi^{(k)}}=-{vol_X (v)}\sum v_i \overline{E_i}.$$
        \item \label{step_continuity} We show that for every $i$, the order of the \textbf{b}-divisor $D_\xi$ along $E_{p_i}$ is $\lim_{k} \ord_{E_{p_i}} D_{\xi^{(k)}}$.
    \end{enumerate}
    By lemma \ref{lem:base-centers} then, we obtain an equality $D_\xi= \lim_{k\to\infty} D_{\xi^{(k)}}$ and the two steps combined finish the proof.
    
    \begin{proof}[Proof of step \ref{step_sequence}]\renewcommand{\qedsymbol}{}
    Since $v$ is not divisorial, the sequence of centers is infinite, and we begin the construction by considering for each $k\in\bbN$ the divisorial valuation
    $$w^{(k)} = \frac{v_1}{\ord_{p_k}(\overline{E_{p_1}})}\cdot \ord_{p_k} .$$
    The normalizing factor $v_1/\ord_{p_k}(\overline{E_{p_1}})$ is chosen so that $w^{(k)}_1=w^{(k)}(E_1)=v_1$ for every $k$. The proximity equalities then imply that, for every $i$, there is a $k_i$ such that $w^{(k)}_i=v_i$ for all $k\ge k_i$.
    Then by Noether's formula, for every $f$ there is a $k_f$ such that $w^{(k)}(f)=v(f)$ for all $k\ge k_f$. Therefore, $w^{(k)}$ converge to $v$. 
    
    The sequence $w^{(k)}$ is in general not nondecreasing, but we will find a suitable subsequence $v^{(k)}=w^{(\ell_k)}$ with this property, which is equivalent to the inequalities
       $v^{(k)}_i \le v^{(k+1)}_i$ for all $i, k \ge 1$. 
    If $p_{k+1}$ is free, the proximity equalities involving $v_i$ and $w^{(k')}_i$ are the same for every $i\le k\le k'$, so $w^{(k)}_i = v_i=w^{(k')}, i=1, \dots, k$ and therefore $w^{(k)}\le w^{(k')}$. 
    Thus, if $v$ has infinitely many free centers, the subsequence formed by those $v^{(k)}=w^{(\ell_k)}$ such that $p_{\ell_k+1}$ is nondecreasing.
    
    Assume now $v$ has only finitely many free centers, and let $p_{\ell_1}$ be the last free center. Then we can find a sequence $\ell_k$ such that $p_{\ell_k+1}, \dots, p_{\ell_{k+1}+1}$ are all the centers $p_i$ proximate to $p_{\ell_k}$ (recall that no center has infinitely many centers proximate to it). 
    The valuation $v$ is Type 5 in the classification of valuations of \cite[Chapter 8]{Cas00}: an irrational valuation associated to an irrational characteristic exponent $\alpha$, and the valuations $v^{(\ell_k)}$ approach $v$ as the approximants of the continued fraction of $\alpha$ approach it: alternatively from above and from below. Then we choose the subsequence of the lower approximants, 
    obtaining a nondecreasing sequence of valuations whose centers are
    contained in the sequence of centers of $v$ and converging to $v$.

	Now let $\ell_k$ be the number of centers of $v^{(k)}$ and consider
	$$
	m_k = \frac{\sum_{i=1}^{\ell_k} v_i^2}{v_{\ell_k}}.
	$$
	The real numbers $m_k$ form a divergent increasing sequence. 
    For any fixed $j$ and for every $k$ we claim the equality
    $\fra_{\xi,m_j}= \fra_{\xi^{(k)},m_j}$ holds. The inclusion
    $\fra_{\xi,m_j}\supset \fra_{\xi^{(k)},m_j}$ holds because $v^{(k)}\ge v$.
    On the other hand, if $f\in \fra_{\xi,m_j}$, let $e_i$ be the multiplicity
    at $p_i$ of the strict transform of $f=0$.
    Noether's equality \ref{lem:noether}
    $$
    v(f)=\sum_{i\ge 1} v_i e_i
    $$
    implies that either $\sum_{i=1}^{\ell_k}v_i e_i\ge m$ (in which case $f\in \fra_{\xi^{(k)},m_j}$) or $e_{\ell_k+1}>0$.
    In the latter case, the proximity inequalities give $e_{\ell_k}\ge 1=v_{\ell_k}/v_{\ell_k}$ 
    and for every $i<\ell_k$, 
    $$
    e_i\ge \sum_{p_{i'}\text{ prox. to }p_i} v_{i'},
    $$
    so by recurrence we get $e_i\ge v_i/v_{\ell_k}$ for all $i\le \ell_k$ and
    therefore  
    $$
	v^{(k)}(f)=\sum_{i= 1}^{\ell_k} v_i e_i \ge \sum_{i= 1}^{\ell_k} v_i^2/v_{\ell_k} =m 
	$$
	i.e., $f\in \fra_{\xi^{(k)},m_j}$ as claimed.
    
    The choice of $m_k$ guarantees that Lemma 
    \ref{lem:am-divisorial} applies and so
    $$
    \length\cO_{X,p}/\fra_{m_k}=\length\cO_{X,p}/\fra_{\xi^{(k)},m_k}= \sum_{i=1}^{\ell_j} \frac{\frac{mv_i}{\sum_{i=1}^{\ell_j} v_i^2}\left(\frac{mv_i}{\sum_{i=1}^{\ell_j} v_i^2}+1\right)}{2}
    $$
    and in particular
    $$
    \vol_X(v)=\frac{1}{\sum_{i\ge 1} v_i^2} = \lim \vol_X(v^{(k)}) \in [0,\infty),
    $$
    so that 
    $$\lim_{k\to \infty}D_{\xi^{(k)}}=-{vol_X (v)}\sum v_i \overline{E_i}$$
    which finishes the proof of step \ref{step_sequence}.
    \end{proof}
    
    \begin{proof}[Proof of step \ref{step_continuity}]\renewcommand{\qedsymbol}{}
    By Proposition \ref{pro:coefficients-dxi}, for every $i$ we have equalities
    \begin{gather*}
        \ord_{E_i}(D_\xi) = -\inf_{f\in \cO_{X,p_1}}\frac{\ord_{E_i}(\pi_i^*(f))}{v(f)},\\
        \ord_{E_i}(D_{\xi^{(k)}}) = -\inf_{f\in \cO_{X,p_1}}\frac{\ord_{E_i}(\pi_i^*(f))}{v^{(k)}(f)}.
    \end{gather*}
    Since $v^{(k)}\le v$ for every $k$, this immediately gives that $\ord_{E_i}(D_\xi)\ge \lim_{k\to\infty} \ord_{E_i}(D_{\xi^{(k)}})$. 
    The opposite weak inequality can be justified directly, for each $i \ge 1$:
    for every $\varepsilon>0$ we will provide an $f\in \cO_{X,p_1}$ and a $k_0$ such that, for every $k\ge k_0$,
    $$
    \frac{ord_{E_i}(\pi_i^*(f))}{v^{(k)}(f)}\le -\ord_{E_i}(D_\xi)+\varepsilon.
    $$
    Indeed, let $f\in \cO_{X,p_1}$ be such that 
    $$
    \frac{ord_{E_i}(\pi_i^*(f))}{v(f)}\le -\ord_{E_i}(D_\xi)+\frac{\varepsilon}{2},
    $$
    and assume (replacing $f$ by a suitable power $f$ if necessary) that $\ord_{E_i}(\pi_i^*(f))\ge 2$.
    Let $k\ge 1$ be such that $v^{(k)}(f)\ge v(f)(1-\frac{\varepsilon}{\ord_{E_i}(\pi_i^*(f))+\varepsilon})$. Then
    $$
    \frac{ord_{E_i}(\pi_i^*(f))}{v^{(k)}(f)}\le \frac{ord_{E_i}(\pi_i^*(f))}{v(f)}
    +\frac{\varepsilon}{\ord_{E_i}(\pi_i^*(f))}\le -\ord_{E_i}(D_\xi)+
    \frac{\varepsilon}{2}+\frac{\varepsilon}{2},
    $$
    as needed.
\end{proof}
     \end{proof}

\begin{corollary}\label{cor:characterize-b-divisorial-surfaces}
If $X$ is a surface, then given a nontrivial norm $\xi\in\frX$, the following are equivalent:
	\begin{enumerate}
		\item \label{part:bdiv} $\xi$ is \textbf{b}-divisorial.
		\item \label{part:nonzerodxi} $D_\xi \ne 0$. 
		\item \label{part:positivevolume} $\red_X(\xi)$ is a curve or $\vol_X(v_\xi)>0$.
	\end{enumerate}
\end{corollary}
\begin{proof}
    We know that \eqref{part:bdiv} implies \eqref{part:nonzerodxi} always,
    and by the preceding Theorem it is clear that \eqref{part:nonzerodxi} implies
    \eqref{part:positivevolume}. Let us prove that \eqref{part:positivevolume}
    implies \eqref{part:bdiv}.
    
    Of course if $\red_X(\xi)$ is a curve then $\xi$ is divisorial, so assume
    $\vol_X(v_\xi)>0$, and keep the notations of the previous proof.
    
    We need to show that if $\div(f)+tD_\xi\ge 0$ then $v(f)\ge t$
    (the converse implication has been proved before).
    We have, for every $k$,
    $$
    \div(f)+t\frac{ \vol_X(v)}{\vol_X(v^{(k)})}D_{\xi^{(k)}} = 
    	\div(f) - t\vol_X(v) \sum_{i=1}^{\ell_k} v_i \overline{E_i} \ge \div(f)+tD_\xi \ge 0,
    $$
    so, because $v^{(k)}$ is divisorial, $v^{(k)}(f)\ge t\vol_X(v)/\vol_X(v^{(k)})$.
    Now, for $k$ large enough,
    $$
    v(f)=v^{(k)}(f) \ge t\frac{ \vol_X(v)}{\vol_X(v^{(k)})}
    $$
    and taking the limit for $k\rightarrow \infty$, we get the desired inequality.
\end{proof}

\begin{remark}
Let us highlight the special properties of valuations on surfaces which are leveraged in the proof of Corollary \ref{cor:characterize-b-divisorial-surfaces} but fail to prove Conjecture \ref{con:characterize-b-divisorial} in higher dimension.

In dimension 2 the sequence of centers of a valuation $v$ is in ordered bijection with the regular local rings dominated by $\cO_v$ and dominating $\cO_{X,\cent_X(v)}$.
In higher dimension, the analogous set of local rings is not totally ordered, and the sequence of centers need not be cofinal with it, and so it does not capture the equivalence class of the valuation (i.e., nonequivalent valuations with the same sequence of centers exist, see \cite{Nun04}).
An attempt to generalize the proof of Corollary \ref{cor:characterize-b-divisorial-surfaces} to prove Conjecture \ref{con:characterize-b-divisorial} would begin by selecting a suitable nondecreasing sequence of divisorial valuations $v^{(k)}$ converging to $v$, and decompositions of the corresponding birational maps in blowups of smooth centers.
If this is feasible, then one would still need to determine coefficients $c(k)$ converging to 1 and such that
$$
c(k)D_{\xi^{(k)}}\ge D_\xi.
$$
This seems totally out of reach. 

In the case of surfaces, $X$-nefness is controlled by the proximity inequalities, which boil down to controlling nefness on the $\bbP^1$'s which arise as exceptional divisors of smooth point blowups.
On the other hand, the values of centers satisfy proximity equalities.
It is the combination of proximity equalities and inequalities that allowed our computation in terms of the volume.

In higher dimension, $X$-nefness becomes a much more challenging question, where a similar approach would require the explicit determination of the nef cone on a sequence of blown up $\bbP^r$-bundles.
This is in general a hard problem, and it is known that in many cases the nef cone is not polyhedral, so the nefness conditions involved in determining the coefficients of $Z(\fra_{\xi,m})$ will not be linear nor related to proximity equalities. 
The sequence of values will continue to satisfy proximity equalities by \ref{lem:noether}, but these will also be harder to control as there is no natural choice for a sequence of centers of the blowups. 

In conclusion, we do not expect an explicit formula like 
$$\frac{\vol_X(v_\xi)}{\vol_X(v^{(k)})} D_{\xi^{(k)}}\ge D_\xi$$
can hold with the needed generality, at least not with a similar proof.
\end{remark}

We end by showing that for every smooth point $p\in X$ on a surfaces, the Seshadri function is lower semicontinuous over the subset of $\frX$ consisting of seminorms centered at $p$. 

\begin{lemma}
    Assume $p_1$ is a smooth point of the surface $X$. Let 
    \[
X_{p_i} \overset{b_{p_i}}{\longrightarrow} 
X_{p_{i-1}}\overset{b_{p_{i-1}}}{\longrightarrow} \dots
\overset{b_{p_1}}{\longrightarrow} X
\]
be a sequence of blowing-ups, write $\pi_i=b_{p_i}\circ b_{p_{i-1}} \circ\dots\circ b_{p_1}$ and let $z=(z_1,z_2)$ be a system of parameters at $p_i\in X_{p_i}$
such that $z_1=0$ is a local equation of the last exceptional divisor $E_{i-1}$.
Then for every nef divisor $D$ on $X$, the map $\varepsilon_i(D,\cdot):\Delta_{\pi_i,p_i,z}\rightarrow\bbR_{\ge 0}\cup \{\infty\}$ is continuous in the interior and on the half-ray $(w_1,0)$ with $w_1>0$.
\end{lemma}
\begin{proof}
    Continuity in the interior is Theorem \ref{thm:toric-seshadri-continuity}.
    Now since $(z_1,z_2)$ is a local system of parameters with $z_1=0$ an equation of $E_{i-1}$,
    $z_2=0$ is transverse to $E_{i-1}$ and so there exist local parameters $u=(u_1,u_2)$ at $p_{i-1}$ such that the map $b_i$ is given locally as $(z_1,z_1z_2)$.
    This implies that the retraction map $\Delta_{\pi_i,p_i,z} \rightarrow\Delta_{\pi_{i-1},p_{i-1},u}$ is an inclusion, given in coordinates as
    $(w_1,w_2)\mapsto(w_1,w_1+w_2)$.
    Since $\varepsilon_{i-1}$ is continuous in the interior of $\Delta_{\pi_{i-1},p_{i-1},u}$, and the ray $(w_1,0)$ is mapped to the interior, the claim follows.
\end{proof}

\begin{theorem}
    \label{thm:continuity-seshadri-surfaces}
    Assume $p$ is a smooth point of the surface $X$. Let $\frX_p=\red^{-1}(p)$.
    The map $\varepsilon(D,\cdot):\frX_p\rightarrow \bbR_{\ge 0}\cup \{\infty\}$ is lower semicontinuous.
\end{theorem}
\begin{proof}
    For every $\xi \in \frX_p$, every model $X_\pi$, point $p\in X_\pi$ and system of parameters $z$, by the previous lemma, the map $\varepsilon(D,\cdot):\Delta_{\pi,p,z}\rightarrow \bbR_{\ge 0}\cup \{\infty\}$ is continuous at $ev_z(\xi)$.
    Then the claim follows by Proposition \ref{pro:seshadri-supremum-models}, as $\varepsilon(D,\cdot)$ is the supremum of continuous functions at a neighborhood of $\xi$.
\end{proof}

{\footnotesize
	\bibliographystyle{plainurl}
	\bibliography{bib_brz}}

\end{document}